\documentclass[12pt]{report}
\usepackage[latin1]{inputenc}
\usepackage{amsmath}
\usepackage{amssymb}
\usepackage{epsfig}
\usepackage{enumerate}
\def\pg{\mathfrak{p}}
\def\Rg{\mathfrak{R}}
\usepackage{amsfonts}
\def\z{\mathbb{Z}}

\newcommand{\ih}{\'\i}  
\newcommand{\eh}{\hspace{.05in}}  
\newcommand{\real}{\mathbb{R}}  
\newcommand{\natl}{\mathbb{N}}  
\newcommand{\C}{\mathbb{C}}
\newcommand{\Z}{\mathbb{Z}} 
 
\newcommand{\CC}{\mathcal{C}} 
\newcommand{\D}{\mathcal{D}} 
\newcommand{\G}{\mathcal{G}} 
\newcommand{\LL}{\mathcal{L}} 
\newcommand{\Q}{\mathcal{Q}}
\newcommand{\R}{\mathcal{R}} 
\newcommand{\T}{\mathcal{T}} 

\newcommand{\ds}{\displaystyle}  
\newcommand{\BE}{\begin{equation}}  
\newcommand{\EE}{\end{equation}}
\newcommand{\m}{\frac{1}{2}}
\newcommand{\tc}{\tilde{c}} 
\newcommand{\tf}{\tilde{f}} 
\newcommand{\tr}{\tilde{r}} 
\newcommand{\tx}{\tilde{x}} 
\newcommand{\deh}{\partial} 
\newcommand{\srho}{\sin\rho}
\newcommand{\af}{\alpha}
\newcommand{\be}{\beta}
\newcommand{\de}{\delta}
\newcommand{\ld}{\lambda}  
\newcommand{\Ld}{\Lambda}  
  
\newcommand{\fns}{\scriptsize} 
\newcommand{\longtr}{\underline{ \ \ \ \ \ \ \ \ \ \ \ \ \ \ \ \ \ \ \ \ \ \ \ }}
\newcommand{\Lim}[1]{\lower5pt\hbox{${{\ds\lim}\atop #1}$}}

\begin{document}
\ \\
{\large{\bf The embedded singly periodic Scherk-Costa surfaces}}   
\\
\\
{\bf Francisco Mart\'\i n $\cdotp$ Val\'erio Ramos Batista}
\\
\\
{\bf Abstract.} We give a positive answer to M. Traizet's open question about the existence of complete embedded minimal surfaces with Scherk-ends without planar geodesics. In the singly periodic case, these examples get close to an extension of Traizet's result concerning asymmetric complete minimal submanifolds of $\real^3$ with finite total curvature.
\\ 
\\
\longtr\longtr\longtr\longtr\longtr
\\
\\ 
{\bf 1. Introduction} 
\\  
\\
The theory of minimal surfaces has been developed for almost two and a half centuries. Along the time, these surfaces acquired several branches of study, classifying them according to their local and global behaviours, total curvature, periodicity, kinds of ends, etc. One of the major interests of study is the class of {\it algebraic} minimal surfaces, conceived as complete minimal surfaces in $\real^3$ obtained from a fundamental piece of finite total curvature by applying a finitely generated translation group $G$ of $\real^3$. The first such surfaces for non-trivial $G$, after the helicoid, were found by H. F. Scherk in 1835 (see [Sk]). It took more than three decades for new other examples to come out, in 1867 by B. Riemann [R] and in 1890 by H. A. Schwarz [Sz].

The presence of Scherk-ends seem to impose strong restrictions to the embedded surface. A characterisation result from W. Meeks and H. Rosenberg in [MR] states that all ends of an embedded doubly periodic algebraic minimal surface {\it must} be of Scherk-type. In the case of a singly periodic algebraic minimal surface, the ends can also be either helicoidal or planar (with finite total curvature). Anyway, for all such surfaces the number of Scherk-ends is {\it even}, four at least, and {\it any} two of them which converge to a common asymptotic plane always contain straight lines. Moreover, to date, the only known singly periodic examples with exactly four Scherk-ends are the members of Scherk's second surface family.
\\
\\
\longtr\longtr\longtr\longtr\longtr
\\
{\footnotesize F. M{\fns ART\'IN}\\
Depto. Geometr\ih a y Topolog\ih a, Universidad de Granada, 18071 \ Granada, Spain (e-mail: \texttt{fmartin@ugr.es}). Mart\ih n's research is partially supported by MCYT-FEDER grant number BMF2001-3489.\\
V. R{\fns AMOS} B{\fns ATISTA}\\
IMECC-Unicamp, P.O.BOX 6065, 13083-970 Campinas, Brazil (e-mail: \texttt{valerio@ime.unicamp.br}).}\\

In fact, it is hard to mention all curiosities about the behaviour of Scherk-ends. Because of that, the study of complete embedded minimal surfaces with Scherk-ends has been split into several particular matters. This paper presents the first examples, called SC surfaces, which give a positive answer to the open problem introduced by M. Traizet in [T1]. It is related to the possibility for existence of complete embedded minimal surfaces in $\real^3$ with Scherk-ends without symmetry curves. In [T2], Traizet gave a positive answer to D. Hoffman and H. Karcher's open problem concerning asymmetric complete minimal submanifolds of $\real^3$ with finite total curvature (see section 5.2 of [HK]).

This paper is devoted to the construction of singly periodic minimal surfaces which are based on the Costa surface in the following sense: one replaces each end of the Costa surface by Scherk-type ends, still keeping the straight lines crossing at the ``Costa-saddle''. This generates a fundamental piece called $S$, which is represented in Figure 1. By successive 180$^\circ$-rotations on the straight lines, one obtains the whole singly periodic SC surface.
\input epsf 
\begin{figure} [ht] 
\centerline{ 
\epsfxsize 8cm  
\epsfbox{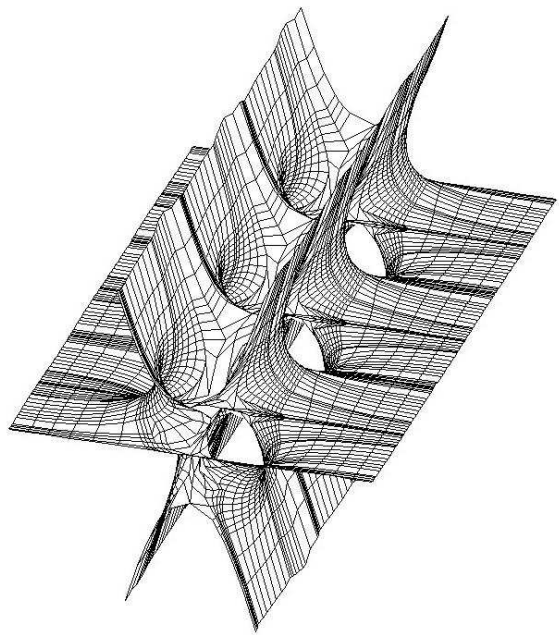}}  
\caption{A singly periodic SC surface.}   
\end{figure} 
 
We try to combine the periodic property of Scherk-ends with the Costa surface and get new examples of minimal surfaces. The main goal of this paper is then to prove the following theorem:
\\
\\ 
{\bf Theorem 1.1.} \it There exists a one-parameter family of complete singly periodic minimal surfaces in $\real^3$ such that, for any member of this family the following holds:
\\
(a) The quotient by its translation group has genus 1.
\\
(b) The whole surface is generated by a fundamental piece, which is a surface in $\real^3$ with border. The fundamental piece has exactly six ends, all of Scherk-type (modulo translation), and a symmetry group generated by a straight line $l_1$ and one straight line segment $l_2$. The line $l_1$ crosses the segment $l_2$ orthogonally in the middle. Both correspond to 180$^\circ$-rotational symmetries of the fundamental piece.
\\
(c) The border of the fundamental piece consists of: (i) two lines parallel to $l_1$, which cross the segment $l_2$ orthogonally at its extremes; and (ii) four isometric curves. These curves provide no extra symmetries for the minimal surface. By successive 180$^\circ$-rotations around the straight lines of the border of the fundamental piece, one generates the whole singly periodic surface.
\\
(d) The SC surfaces are embedded in $\real^3$.\rm
\\

Sections 8 and 9 of this paper show that the period problems can be solved, while (d) is proved in Section 10. The existence proof of the SC surfaces began during the second author's doctoral studies in Germany (see [RB1]), supported by CAPES and DAAD, {\it under} the supervision of Prof. Dr. H. Karcher. We thank him for his valuable advisements, which made possible the realisation of this present work. The embeddedness of the SC surfaces turned out to be a second research project, this time supported by FAPESP grant numbers 00/07090-5, 02/00865-7 and 03/01398-6. 
\\
\\
{\bf 2. Background}
\\
\\
In this section we state some well known theorems on minimal surfaces. For details we refer the reader to [K], [LM], [N] and [O]. In this paper all surfaces are supposed to be regular.   
\\
\\
{\bf Theorem 2.1.} (Weierstrass representation). \it Let $S$ be a minimal surface in $\real^3$ and $R$ the underlying Riemann surface of $S$. Let $dh$ be a meromorphic 1-form on $R$ and $g:R\to\hat\C:=\C\cup{\infty}$ a meromorphic function. Then $X:R\to\real^3$ given by 
\[
   X(p):=Re\int^p(\phi_1,\phi_2,\phi_3),\eh\eh where\eh\eh
   (\phi_1,\phi_2,\phi_3):=\m\biggl(g-g^{-1},i(g+g^{-1}),2\biggl)dh,
\]
is a conformal regular minimal immersion provided the poles and zeros of $g$ coincide with the zeros of $dh$. Conversely, every regular conformal minimal immersion $X:R\to\real^3$ can be expressed in this form for some meromorphic function $g$ and meromorphic 1-form $dh$.\rm
\\
\\
{\bf Definition 2.1.} The pair $(g,dh)$ is the \it Weierstrass data \rm on $R$ of the minimal immersion $X:R\to X(R)=S\subset\real^3$.
\\
\\
{\bf Definition 2.2.} A complete, orientable minimal surface $S$ is \it algebraic \rm if it admits a Weierstrass representation such that $R=\bar{R}\setminus\{p_1,\dots,p_r\}$, were $\bar R$ is compact, and both $g$ and $dh$ extend meromorphically to $\bar R$. 
\\
\\
{\bf Definition 2.3.} An \it end \rm of $S$ is the image of a punctured neighbourhood $V_p$ of a point $p\in\{p_1,\dots,p_r\}$ such that $(\{p_1,\dots,p_r\}\setminus\{p\})\cap\bar{V}_p=\emptyset$. The end is \it embedded \rm if this image is embedded for a sufficiently small neighbourhood of $p$. 
\\
\\
{\bf Theorem 2.2.} \it Let $S$ be a complete minimal surface in $\real^3$. Then $S$ is algebraic if, and only if, it can be obtained from a piece $\tilde{S}$ of finite total curvature by applying a finitely generated translation group $G$ of $\real^3$.\rm 
\\

From now on we consider only algebraic surfaces. The function $g$ is the stereographic projection of the Gau\ss \ map $N:R\to S^2$ of the minimal immersion $X$. This minimal immersion is well defined in $\real^3/G$, but is allowed to be a multivalued function in $\real^3$. The function $g$ is a covering map of $\hat\C$ and the total curvature of $\tilde{S}$ is $-4\pi$deg$(g)$. 
\\
\\
{\bf 3. The Costa surface}
\\
\\
We describe the Costa surface, which is the starting point of our constructions. Details are found in [C], [HM] and [K].
\\
\\
{\bf Theorem 3.1.} (The Costa surface). \it Let $\bar R$ be the square torus of which the algebraic equation is  
\[
   \wp'^2=i\wp(\wp-1)(\wp+1).
\] 
For some positive $\mu$ define $g=\mu\wp'$ and $dh=d\wp/(\wp^2-1)$. Then there exists a unique positive $\mu_0$ such that, for $\mu=\mu_0$, $(g,dh)$ is the Weierstrass data on $R=\bar{R}\setminus\wp^{-1}(\{-1,1,\infty\})$ of a complete minimal embedding of $R$ in $\real^3$.\rm    
\\

Figure 2(a) represents the image of the minimal embedding referred to by the previous theorem.
\\

\input epsf
\begin{figure} [ht]
\centerline{
\hspace{-1in}
\epsfxsize 9cm
\epsfbox{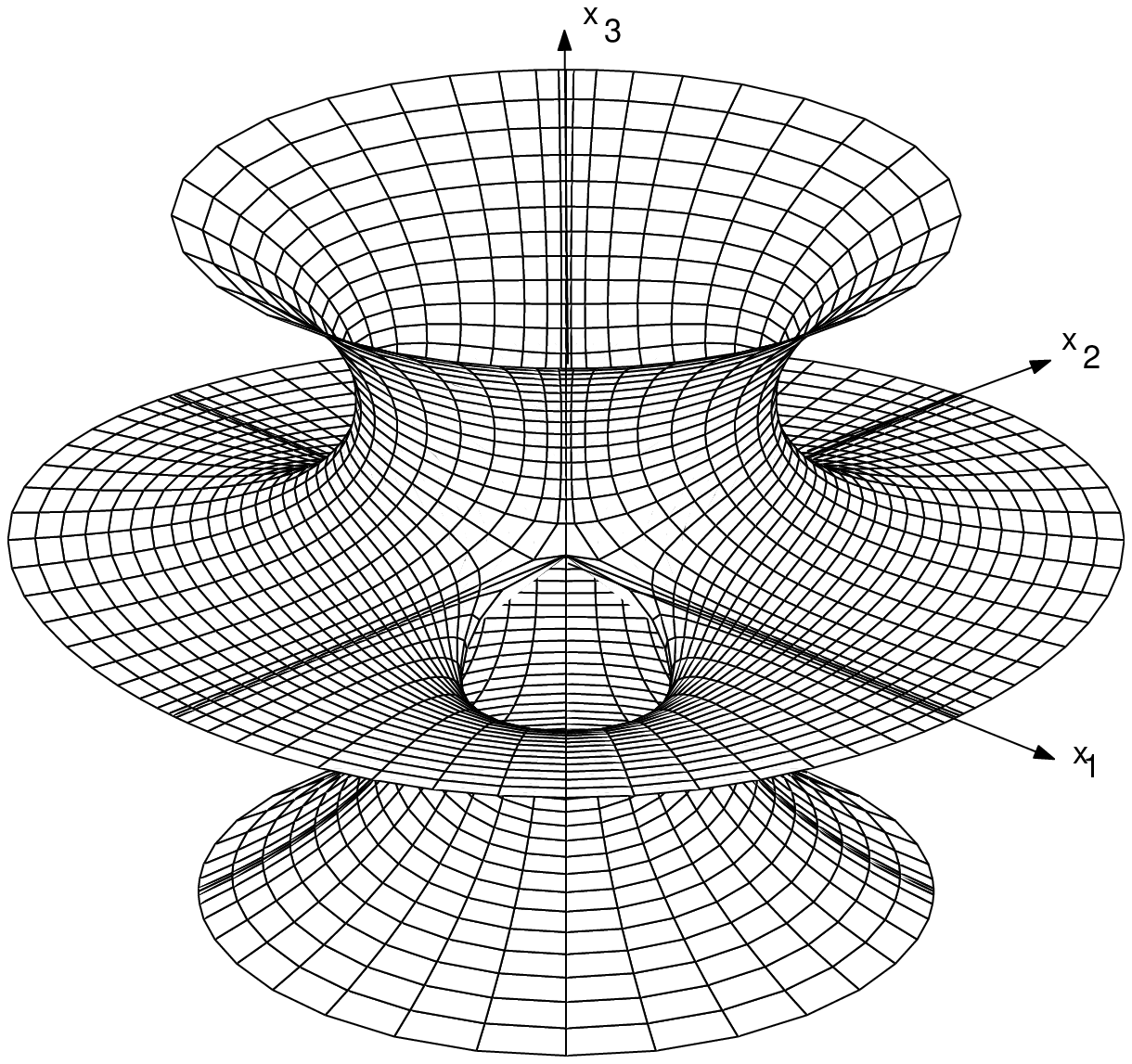}
\epsfxsize 8cm
\hspace{-.1in}
\epsfbox{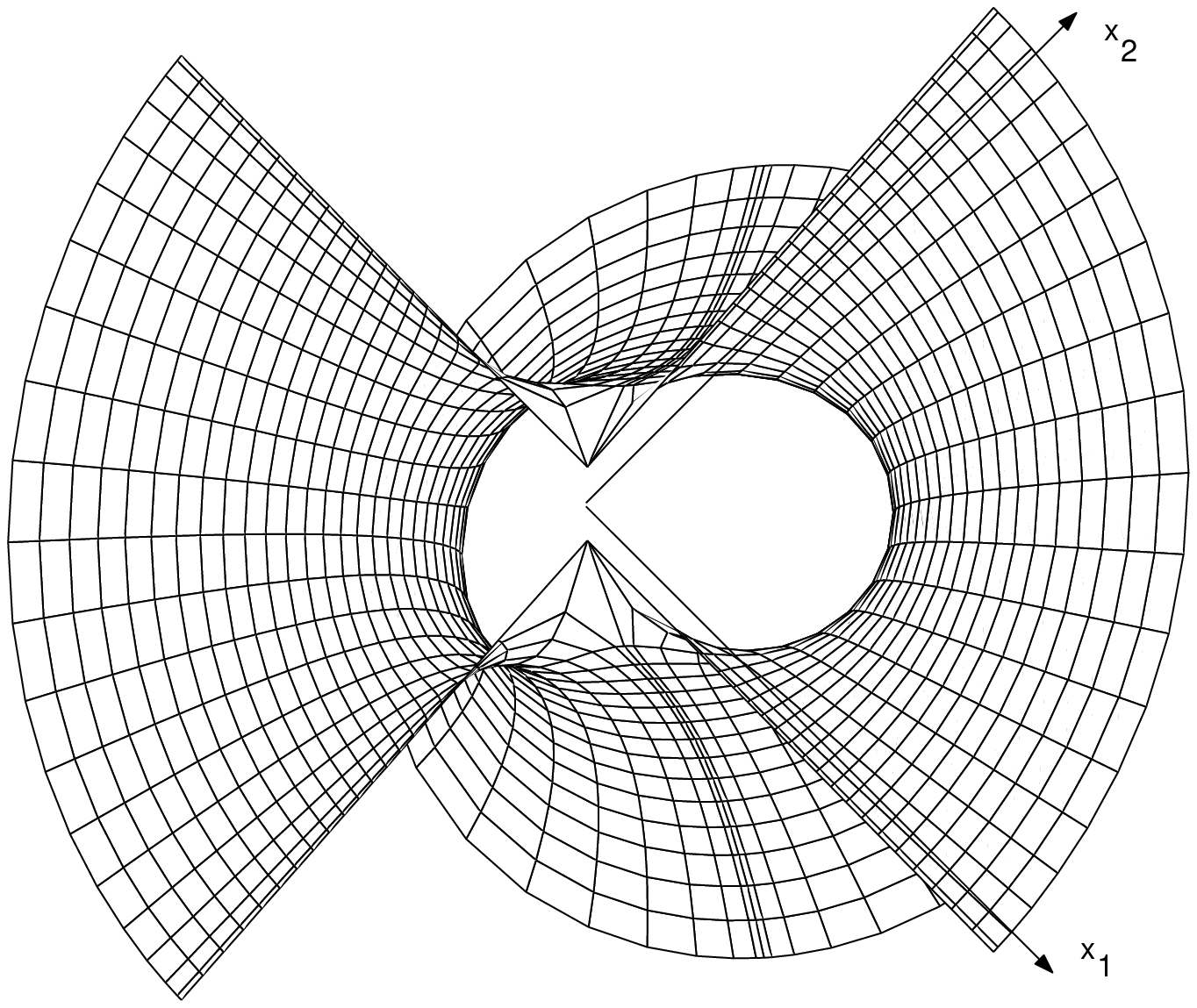}}
\hspace{1.6in}(a)\hspace{3in}(b) 
\caption{(a) The Costa surface in $\real^3$; (b) the case $\mu>\mu_0$.}
\end{figure}

After a suitable rigid motion in $\real^3$ we can position the Costa surface in such a way that the generators of its symmetry group are:
\\
\\
\begin{tabular}{llc} 
 $ \hspace{1.7in} $&$\sigma_1:=(x_1,x_2,x_3)\to(-x_2,-x_1,x_3);$\\
 $                $&$\sigma_2:=(x_1,x_2,x_3)\to(x_2,x_1,x_3);$\\
 $                $&$\sigma_3:=(x_1,x_2,x_3)\to(-x_1,x_2,-x_3)\eh\eh{\rm and}$\\
 $                $&$\sigma_4:=(x_1,x_2,x_3)\to(x_1,-x_2,-x_3).$  
\end{tabular} 
\ \\
\\

Notice that $\sigma_2=\sigma_3\circ\sigma_1\circ\sigma_3$ and $\sigma_4=\sigma_1\circ\sigma_3\circ\sigma_1$. We call $\tilde{G}$ the group of symmetries of the Costa surface. In our case,
\[
   \tilde{G}=<\sigma_1,\sigma_3>.
\]

We remark that the Costa surface is invariant under a 180$^\circ$-rotation around the $x_3$-axis. This rotation can be given by $\sigma_1\circ\sigma_2$. Now consider the shaded region of the punched square torus $R$ represented in Figure 3. This region can be thought of as a punched square $\Q$ in the plane. Regarding the referred to constant $\mu_0$ in Theorem 3.1, the same Weierstrass data still provide minimal immersions (with boundary) of $\Q$ in $\real^3$, for $\mu>\mu_0$. These are also embeddings, and the proof of this fact follows the same arguments as in [K] for the embeddedness proof of the Costa surface. Figure 2(b) illustrates an ``open half-Costa surface'', which can still be positioned in $\real^3$ in order to have the symmetries $\sigma_1$ and $\sigma_2$.
\\
\input epsf
\begin{figure} [ht]
\centerline{
\epsfxsize 9cm
\epsfbox{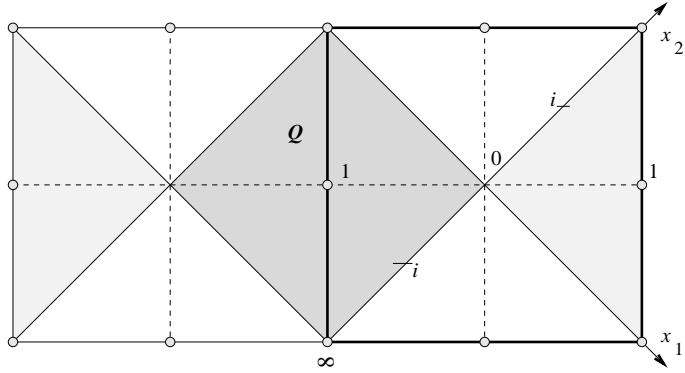}}
\caption{Two copies of the square torus with some values of $\wp$.}
\end{figure}
\ \\
{\bf 4. The symmetries of the SC surfaces and the elliptic $z$-functions}
\\ 
\\  
Let us consider the surface represented in Figure 4 and call it $S$. The quotient by its translation group, followed by a compactfication of its ends, generates a torus $T$ (see Figure 5(a)). The straight lines on Figure 4 are 180$^\circ$-rotational symmetries of $S$. We consider that two lines of $S$ coincide with the $x_1$- and $x_2$- axis, respectively. The others will be parallel to $x_1$. 
\input epsf 
\begin{figure} [ht] 
\centerline{ 
\epsfxsize 7cm  
\epsfbox{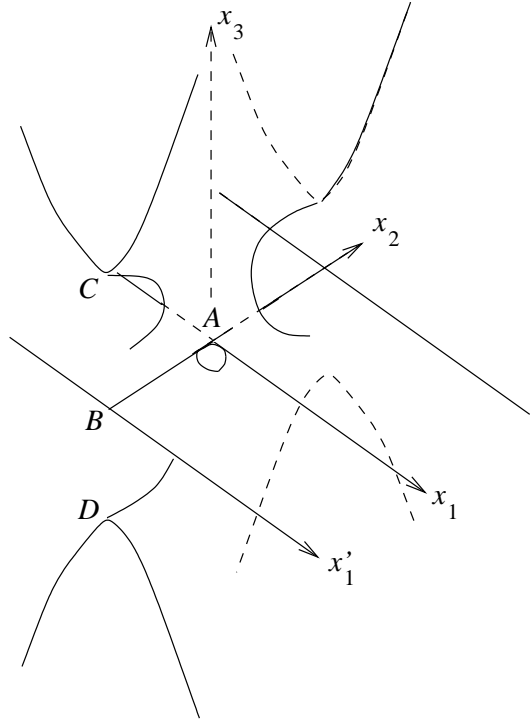}}  
\caption{The surface $S$ with special points on it.} 
\end{figure} 
 
These rotational symmetries turn out to be involutions of $T$. For each one of the straight lines, the fixed point set of the 180$^\circ$-rotation around it has just one component, namely the line itself. Because of that, $T$ must be a rhombic torus.  
 
On the surface $S$, the axes $x_1$ and $x_2$ cross at a point $A$. The axes $x_{1}^{\prime}$ and $x_2$ cross at a point $B$. The surface $S$ has another 180$^\circ$-rotational symmetry, which is around the $x_3$-axis, and the quotient of $T$ by this symmetry is $S^2$. After we fix an identification of $S^2$ with $\hat\C$, this defines an elliptic function $z:T\to S^2$. 
 
Let us call $\T$ the torus $T$ punctured at some points to be defined later. Suppose that the surface $S$ is a minimal immersion of $\T$ in $\real^3$. Then, the unitary normal vector on $S$ will be vertical at $A$ and $B$. Moreover, there will be two other points at which the normal vector will be vertical as well. These we call $C$ and $D$. 
 
Now consider Figure 5(a) and the points of the torus $T$ represented there. These points correspond to special points of $S$ (they were given the same names). Let $z:T\to S^2$ be the elliptic function with $z(A)=0$, $z(B)=\infty$ and $z(C)=e^{i\rho}$, for some $\rho\in(-\pi/2,\pi/2)$. The most important values of $z$ on $T$ are represented in Figure 5(b). 
\input epsf 
\begin{figure} [ht] 
\centerline{ 
\epsfxsize 16cm  
\epsfbox{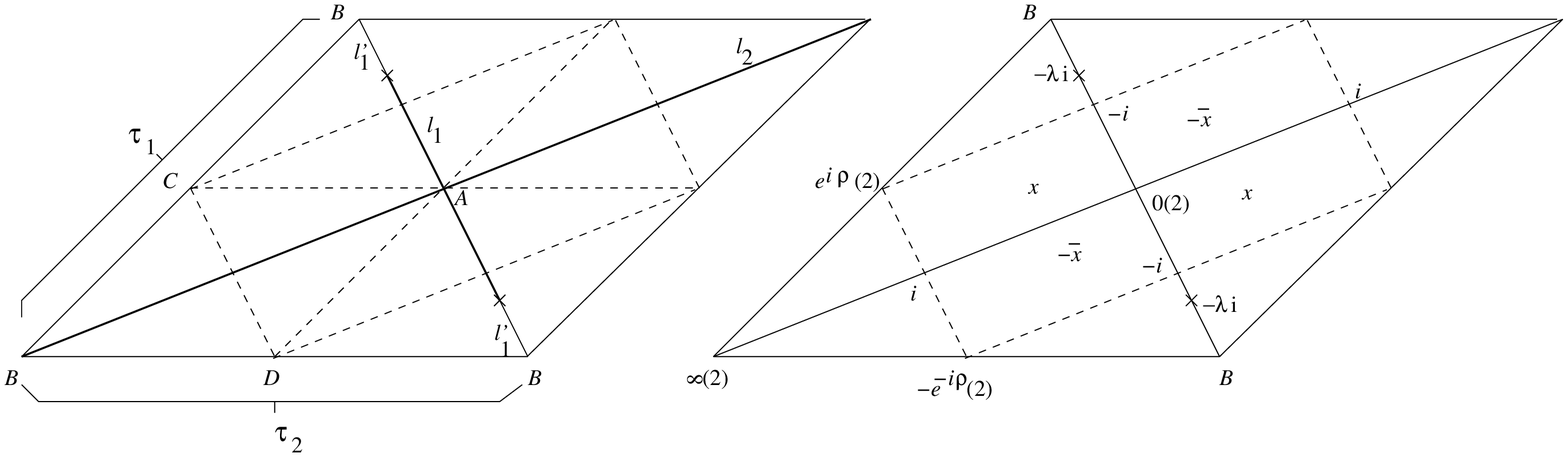}}  
\hspace{1in}(a)\hspace{2.5in}(b) 
\caption{(a) The torus $T$; (b) special values of $z$ on $T$.} 
\end{figure} 
 
On Figure 5(a), the bold lines $l_1$ and $l_2$ represent the $x_1$- and $x_{2}$-axes, respectively (compare it with Figure 4). The $x_{1}^{\prime}$ axis from Figure 4 is indicated by $l_{1}^{\prime}$ in Figure 5(a). On this figure, notice that the points marked with a $\times$ correspond to Scherk-ends of $S$, where $z$ assumes the value $-\ld i,\ld\in(0,\infty)$. The other Scherk-ends correspond to $z=x$ and $z=-\bar{x}$, represented in Figure 5(b). Therefore, we define $\T:=T\setminus z^{-1}(\{x,-\bar{x},-i\ld\})$.
 
Now we are going to summarise important properties of the function $z$ (see Figures 5(a) and 5(b)). The function $z$ is pure imaginary exactly on $l_1,l_{1}^{\prime}$ and $l_2$. On $l_2$, the imaginary part of $z$ is always positive. On $l_1$ and $l_{1}^{\prime}$, the imaginary part of $z$ is always negative. On the dashed lines of Figure 5(b), and nowhere else, $z$ is unitary ($|z|\equiv 1$). 
 
From now on we shall consider $T$ as a general member of the family of tori, of which the algebraic equations are 
\BE            
   z'^2=-iz(z-e^{i\rho})(z+e^{-i\rho}),\eh\eh\rho\in(-\pi/2,\pi/2). 
\EE 
We label $\tau_1$ and $\tau_2$ as two generators of the first homology group of $T$, denoted  $H_1(T,\z).$ \\

\noindent {\bf 5. The Weierstrass data $(g,dh)$ in terms of $z$} 
\\
\\ 
We are supposing that the surface $S$ is a minimal immersion of $\T$ in $\real^3$. In Figure 5(a), the point $D$ corresponds to the image of $C\in S$ under the 180$^\circ$-rotation around the $x_2$-axis (see Figure 4). On the surface $S$, the unitary normal vector at $A,B,C$ and $D$ is vertical. We consider that $g(B)=\infty$ and $g=0$ at $A,C$ and $D$. At the middle Scherk-ends of $S$, marked with $\times$ on Figure 5(a), $g$ takes the value $\infty$ as well. These ends correspond to some value $z=-\ld i$, where $\ld$ is a positive real. We do not expect the unitary normal vector on $S$ to be vertical at any other points except the ones mentioned here. Because of that, one obtains the following relation between $g$ and $z$, for a certain positive constant $c$: 
\BE 
   g=cw,\eh\eh{\rm where}\eh\eh w=\frac{z'}{z+\ld i}. 
\EE  

If we suppose that the surface $S$ is minimal, then the corresponding differential $dh$ for the Weierstrass data $(g,dh)$ must have a single zero at every point of $S$ where $g$ is vertical, {\it except} at the two Scherk-ends corresponding to $z=-i\ld$. These ends are the ones which come from the planar end of the Costa surface. 
 
Two of the other Scherk-ends of $S$ must correspond to some complex value $z=x$, distinct from the special $z$-values we have already mentioned. This means, $x\not\in\{0,\infty,e^{i\rho},-e^{-i\rho},-\ld i\}$. Moreover, we must have $x\not\in i\real$ as well, since $i\real$ is the image of the straight lines of $S$ under $z$. By the way, exactly because of these straight lines, the two remaining Scherk-ends of $S$ must be at $z=-\bar{x}$. These arguments lead to the following conclusion:  
\BE 
   dh=\frac{dz}{(z-x)(z+\bar{x})}=\frac{dz}{(z-ib)^2-a^2},\eh\eh 
   a+ib:=x,\{a,b\}\subset\real. 
\EE 
 
Of course, both sides of (3) are a priori just proportional, but on the straight lines of $S$, where $z(t)=it,t\in\real$, we want $Re\int dh$ to be zero. This is already provided by (3). 
 
Now we have concrete Weierstrass data $(g,dh)$ on the torus $T$. From (1), (2) and (3), it is easy to show that $\frac{\ds dg}{\ds g}\cdot dh$ is pure imaginary on $z(t)=it,t\in\real$, which proves that $S$ really has the straight lines we supposed at the beginning. 
\\ 
\\ 
{\bf 6. Residue problem: conditions on the variables $c,x$ and $\ld$}  
\\
\\ 
We are supposing that $S$ is a minimal immersion of $\T$ in $\real^3$. In the previous section, we obtained its corresponding Weierstrass data (2) and (3) in terms of the symmetries that $S$ was assumed to have. Therefore, (2) and (3) give us a minimal immersion in $\real^3$ of the torus $T$ punctured at $z^{-1}(\{-\ld i,x,-\bar{x}\})$. The image of this immersion in $\real^3$ has the same symmetries and Scherk-ends as the surface $S$ indicated by Figure 1. 
 
Nevertheless, we must introduce some conditions involving the variables $c,x$ and $\ld$, for a given $\rho\in(-\pi/2,\pi/2)$. This is because the minimal immersion must lead to an image in $\real^3$ with no other periods, except a {\it single one} in the $x_2$-direction (see Figure 1). Therefore, all Scherk-ends must have the same length, modulo translation. Moreover, by looking at Figure 4 we realise that the stereographic projection of the unitary normal vector at the ends $z=x$ and $z=-\bar{x}$ {\it must} be real. Otherwise $S$ will have periods in the $x_1$- or $x_3$-directions. 
 
Before we calculate the residues, let us establish {\it necessary} conditions for $g|_{z=x}$ to be real. For example, $g|_{z=x}$ applied to (2) must give us 
\BE 
   \frac{-ix(x-e^{i\rho})(x+e^{-i\rho})}{(x+\ld i)^2}=r^2,r\in\real_+^*. 
\EE 
 
Now, from (1) and (2), notice that $z\to-\bar{z}$ implies $g\to\pm\bar{g}$. Therefore, $g|_{z=-\bar{x}}$ applied to (2) leads to the same condition (4). By using the Cardano-Tartaglia formula we get 
\[
   x=\sqrt{3}(u-v)/2-i[(u+v)/2+R/3],\eh\eh{\rm where}
\]
\[
   R=-2\srho-r^2,\sigma=1-2\ld r^2,\tau=-\ld^2 r^2, p=\sigma-\frac{R^2}{3}, 
   q=\frac{2R^3}{27}-\frac{\sigma R}{3}+\tau, 
\]
\BE
   \D=\frac{p^3}{27}+\frac{q^2}{4}, u=\sqrt[3]{-\frac{q}{2}+\sqrt{\D}}
   \eh\eh{\rm and}\eh\eh v=-\frac{p}{3u}.
\EE

Of course, we do not always have $u\in\real$, but for certain ranges of $\rho,\ld$ and $r$ it will follow that $\D\ge 0$. From (4) we easily see that one of the roots is pure imaginary (with positive imaginary part), and the other two are interchangeable by $z\to -\bar{z}$. Because of this, without loss of generality we fix $a=Re\{x\}>0$.
\\
\\
1) Residue at $z=-\ld i$.  
 
Take a positive $\de$ and the curve $z(t)=-\ld i+\de e^{it},t\in[0,2\pi]$. We recall that $x=a+ib$, for reals $a$ and $b$. An easy calculation shows that, up to the sign of the square root, we have: 
\\
\[ 
   \lim_{\de\to 0}Re\int_{z(t)}\phi_1=0, 
\] 
\BE 
   \lim_{\de\to 0}Re\int_{z(t)}\phi_2=\frac{c\pi}{(\ld+b)^2+a^2} 
                                         [\ld(\ld^2+1+2\ld\srho)]^{1/2}, 
\EE  
\[ 
   \lim_{\de\to 0}Re\int_{z(t)}\phi_3=0. 
\] 
\ \\ 
2) Residues at $z=x$ and $z=-\bar{x}$. 
 
We have just seen that the absolute value of $r$ in (4) remains invariant under the change $x\to-\bar{x}$. Moreover, from (3) the same holds for $dh/dz$. Because of this, the residues at $z=x$ and $z=-\bar{x}$ will have the same absolute values. This is a fact we have already expected since the straight lines on $S$ imply that the Scherk-ends at $z=x$ and $z=-\bar{x}$ must be congruent. Thus, it is sufficient to analyse the residue at $z=x$. 
 
Take a positive $\de$ and the curve $z(t)=x+\de e^{it},t\in[0,2\pi]$. An easy calculation shows that, up to the sign of the square root, we have: 
\[ 
   \lim_{\de\to 0}Re\int_{z(t)}\phi_1=0, 
\] 
\BE 
   \lim_{\de\to 0}Re\int_{z(t)}\phi_2= 
   \frac{\pi}{2a}\biggl(cr+\frac{1}{cr}\biggl), 
\EE 
\[ 
   \lim_{\de\to 0}Re\int_{z(t)}\phi_3=0. 
\] 
 
Now we are going to make (6) equal to (7) but first, to simplify the formulae we are going to derive, let us define for {\it positive} t:
\BE 
   f(t):=[t(t^2+1+2t\srho)]^{1/2}\eh\eh{\rm and}\eh\eh
   \tf(t):=[t(t^2+1-2t\srho)]^{1/2}. 
\EE 
 
Thus, all lengths of the Scherk-ends (modulo translation) will be equal if, and only if
\BE \label{eq:macacoloco}
   c^2=\frac{(\ld+b)^2+a^2}{2af(\ld)r-r^2[(\ld+b)^2+a^2]}. 
\EE 
 
We should be careful at this point and analyse the conditions for $\rho,\ld$ and $r$ which make positive the denominator of (9). However, in the next section we shall get another expression for $c^2$, which solves the period problems for $H_1(T,\z)$. This expression, in (14), will be positive for certain values of $\rho,\ld$ and $r$, and continuous for these variables. 

A solution for the period and residue problems will then exist if, and only if both (9) and (14) simultaneously hold, for certain values of $\rho,\ld$ and $r$. For these values (9) will be automatically positive. Therefore, we skip the analysis of the conditions which lead to a positive denominator in (9).
\\
\\
{\bf 7. The period problems} 
\\
\\ 
We have already analysed the periods at the punctures $z=-\ld i,z=x$ and $z=-\bar{x}$ of the torus $T$ in the last section. Therefore, we need to verify the remaining periods on curves of $H_1({\cal T} ,\z)$. 
Recall that we have the following Weierstrass data:
$$g=c \, w, \quad \mbox{where} \; \; w=\frac{z'}{z+i \ld}; \; \ld \in \real_+,$$
$$dh=\frac{dz}{(z-x)(z+\overline{x})},$$
defined on the tori $T$ with algebraic equation:
$$z'^2=-i z (z-e^{i \rho})(z+e^{-i \rho}), \quad \rho \in \left(-\frac{\pi}{2},\frac{\pi}{2}\right).$$
The ends correspond to the points in $z^{-1}( \{-i \ld, x,-\overline x\})$, that we label $-i \ld_+$, $-i \ld_-$, $x_+$, $x_-$, $-\overline x_+$, and $-\overline x_-$. Recall that ${\cal T}=T-\{-i \ld_+, -i \ld_-, x_+, x_-, -\overline x_+, -\overline x_-\}.$
\input epsf 
\begin{figure} [ht] 
\centerline{ 
\epsfxsize 8cm  
\epsfbox{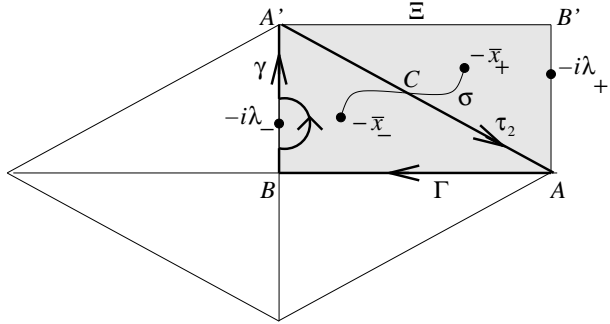}}  
\caption{The rectangle $\Xi$ and the curve $\sigma$.}
\label{fig:FIGURE6B}
\end{figure} 

\begin{figure} [ht] 
\centerline{ 
\epsfxsize 10cm  
\epsfbox{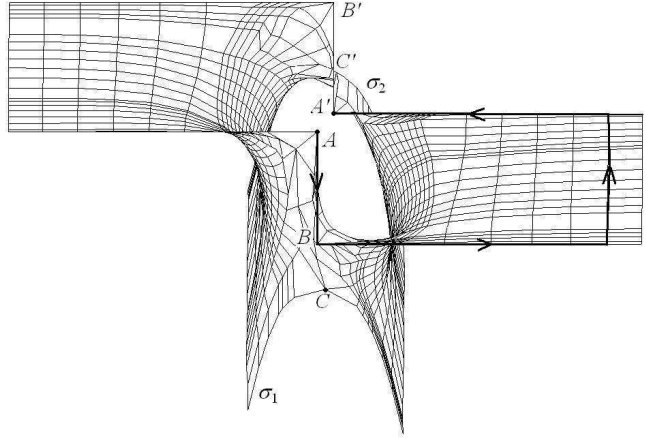}}
\caption{The image of  ${\cal A}:=\Xi- \sigma$ under the minimal immersion $X$. Observe that $\sigma$ is lifted to two curves, $\sigma_1$ and $\sigma_2$, which differ by the period $\vec p$.}  
\label{fig:FIGURE6C}
\end{figure}

If we consider the rectangle $\Xi$ in the plane given by Figure \ref{fig:FIGURE6B}, and we remove a simple curve $\sigma$, joining $-\overline{x}_+$ and $-\overline{x}_-$ and passing through $C$, then the remaining domain has the conformal type of an annulus (with one boundary) and twice punctured at the boundary. The image of  ${\cal A}:=\Xi- \sigma$ under the minimal immersion $X=Re  \int\left(\phi_1,\phi_2,\phi_3 \right)$ has the shape exhibited in Figure \ref{fig:FIGURE6C}.
\input epsf 

Notice that $X$ is well defined on $\cal A$ because the integrals around $x_+$ and $x_-$ have opposite signs. Moreover, observe that $\sigma$ is lifted to two curves, $\sigma_1$ and $\sigma_2$, which differ by the period
$$\vec p=\left(0,\frac\pi{2a} \left(c\, r+\frac{1}{c \, r} \right),0 \right), \quad \mbox{where $c$ is given by (\ref{eq:macacoloco}).}$$
Furthermore, the distance between the two parallel half lines of either the end $-i \ld_+$ or $-i \ld_-$ is $\frac12 \|\vec p\|.$

\begin{figure}[hbtp]
	\begin{center}
		\includegraphics[width=0.60\textwidth]{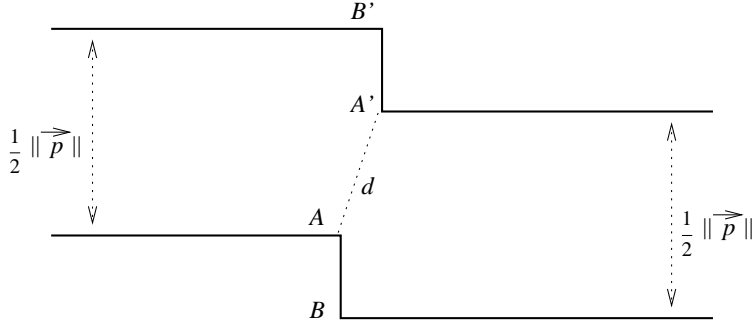}
	\end{center}
	\caption{The intersection of $X({\cal A})$ with the plane $x_3=0$.}
	\label{fig:FIGURE6D}
\end{figure}

If we consider the intersection of $X({\cal A})$ with the plane $x_3=0$, then we see a configuration of half lines and segments like in Figure \ref{fig:FIGURE6D}. A necessary and sufficient condition to close the period problem of $X$ is that {\itshape the distance $d$ vanishes.} In other words, if $d=0$ then we can extend our immersion by Schwarz reflexion principle in order to obtain a complete minimal surface properly immersed in $\real^3$. This complete minimal surface is invariant under the horizontal translation by vector $\vec p.$ Hence, labelling $\bf P$ the group of translations generated by $\vec p$, it is clear that $X$ induces a complete minimal immersion ${\cal X}: {\cal T} \longrightarrow \real^3/{\bf P}.$
At this point, we consider the curves $\gamma(t)$ and $\Gamma(t)$ defined as follows. Given a positive $\de<\ld$, we can write down explicit equations for $z(\Gamma)$ and $z(\gamma)$: 
\BE 
   z(\Gamma(t))=it,0<t<\infty,\eh\eh{\rm and} 
\EE 
\BE  
  z(\gamma(t))=\begin{cases}  
 \gamma_1(t)=-it,\ld+\de<t<\infty, \cr  
   \gamma_\de(t)=-\ld i+\de ie^{it},0<t<\pi,\cr  
   \gamma_2(t)=-it,0<t<\ld-\de.\cr \end{cases} 
\EE  
 Observe that, in the homology group of our punctured torus ${\cal T}$, we have:
$$  (\Gamma \cup \gamma)+\tau_2=\alpha, $$
where $\tau_2$ is one of the generators of the homology of the compact torus and $\alpha$ is a cycle around the end $-\overline{x}_-.$ Taking the previous arguments into account, the periods of our surface are closed if, and only if
\BE \label{eq:cabrito} Re \int_{\Gamma \cup \gamma} (\phi_1,\phi_2)=(0,0). \EE
Notice that (\ref{eq:cabrito}) is equivalent to $d=0.$ Indeed, if $X(\Gamma \cup \gamma)$ is a closed curve in $\real^3$ this means that  $A$ and $A'$ are mapped in the same point of space.
\\
\\
\underline{Remark 1}: the points $B$ and $C$ lie on the same vertical. Indeed, we have two vertical lines of symmetry: one of them lies on the vertical of $A$ and the other one lies on the vertical of $B$. The composition of the $180^o$ rotation around both straight lines gives us the period translation, so they induce the same transformation on $\cal T$. Furthermore, this holomorphic involution corresponds to the elliptic involution $(z,z')  \mapsto (z,-z').$ The involution fixes the points $A$, $B$ and $C$. Thus $C$ is on the same vertical of either $A$ or $B$. The first option is absurd, then $B$ and $C$ are on the same vertical.
\\

At a point $(c,\ld,x)$, a necessary and sufficient condition for (\ref{eq:cabrito}) to hold is $\overline{\int_{\Gamma\cup\gamma}gdh}=\int_{\Gamma\cup\gamma}\frac{\ds dh}{\ds g}$. From (2) and (12), we rewrite this condition as 
\BE
   c^2=\frac{\int_{\Gamma\cup\gamma}w^{-1}dh}{\overline{\int_{\Gamma\cup\gamma}wdh}}.
\EE
 
Now we shall use explicit integrals in (13). The straight lines of $S$ are represented by $z(t)=it,t\in\real$. Therefore, from (2), (3) and (8) one rewrites (13) as follows: 
\BE
   c^2=\frac{\ds\int_{0}^{\infty}\frac{(t-\ld)dt}{f(t)[(t+b)^2+a^2]}
          +i\int_{0}^{\infty}\frac{(t+\ld)dt}{\tf(t)[(t-b)^2+a^2]}}
            {\ds\int_{0}^{\infty}\frac{f(t)dt}{(\ld-t)[(t+b)^2+a^2]}
          +i\frac{\pi f(\ld)}{(\ld+b)^2+a^2}-i\int_{0}^{\infty}\frac{\tf(t)dt}{(t+\ld)[(t-b)^2+a^2]}}.
\EE

Notice the term $(\ld-t)^{-1}$ in one of the integrands of (14). The integral is even though finite and equals to its Cauchy principal value. Now we can split (14) in the following {\it two} conditions
\BE
   c^2=\frac{\ds\int_{0}^{\infty}\frac{(\ld-t)dt}{f(t)[(t+b)^2+a^2]}}
            {\ds\int_{0}^{\infty}\frac{f(t)dt}{(t-\ld)[(t+b)^2+a^2]}}\eh\eh{\rm and}\eh\eh
   c^2=\frac{\ds\int_{0}^{\infty}\frac{(t+\ld)dt}{\tf(t)[(t-b)^2+a^2]}}
         {\ds\frac{\pi f(\ld)}{(\ld+b)^2+a^2}-\int_{0}^{\infty}\frac{\tf(t)dt}{(t+\ld)[(t-b)^2+a^2]}}.
\EE

Together with (9), we now have {\it three} distinct equations for $c^2$, which depend only on $r,\rho$ and $\ld$. Figure 9 shows numeric simulations of these three equations for $\srho=-0.23$ and $r$ in the interval $[0.65,0.75]$. Therefore, we have computational evidence that the period problems can be solved. A formal proof of this fact will be shown in the next sections. For now, we give here a summary of this proof:  
\\
(a) We first define $c_1:=c^2(\rho,\ld,r)$, where $c^2$ is given in (9), $c_2:=c^2(\rho,\ld,r)$, where $c^2$ is the first condition at (15) and finally $c_3:=c^2(\rho,\ld,r)$, where $c^2$ is the second. At a certain point, we shall replace $c_3$ by another equivalent L\'opez-Ros parameter $\tc_3$. The reason for it will be explained later in details.      
\\
(b) Afterwards we show that for negative $\rho\cong 0$ the functions $c_1$ and $\tc_3$ can be estimated by two functions each. In other words, $c^{-}_{1}<c_1<c^{+}_{1}$ and $\tc^{-}_{3}<\tc_3<\tc^{+}_{3}$, and these inequalities are valid for $\R:=\ld r\in(0,r_0)$, $r_0$ to be specified later. Moreover, the graphs of each pair of these estimating functions intersect at most at a single point in $(0,r_0)$. This fact is illustrated by Figure 10.
\\   
(c) We shall conclude that $c_2<\min\{c^{-}_{1},\tc^{-}_{3}\}$ in $(0,r_0)$ for $\ld=1$ and $c_2>\max\{c^{+}_{1},\tc^{+}_{3}\}$ for $\ld=\infty$ in $(0,r_0)$. This will finally prove the existence of a point $(\ld^*,\R^*)\in(1,\infty)\times(0,r_0)$ at which $c_1=c_2=c_3$. 

This strategy is restricted to the case of negative $\rho\cong 0$. Otherwise the $\ld$- and $r$-intervals can drastically change. This really happens according to some numeric computations, and they also indicate solutions for several values of positive and negative $\rho$. In any case, the Cardano-Tartaglia formulae for $x$ in (4) get more complicated for $\rho\ne 0$ and it becomes very difficult to handle the derivatives $\frac{\ds\deh c_j}{\ds\deh r},\eh j=1,2,3$. Because of that, we restrict our theoretical solution of the period problems to the proof of the following proposition, which will be shown in Section 9:

\input epsf  
\begin{figure} [ht]  
\centerline{  
\epsfxsize 7cm   
\epsfbox{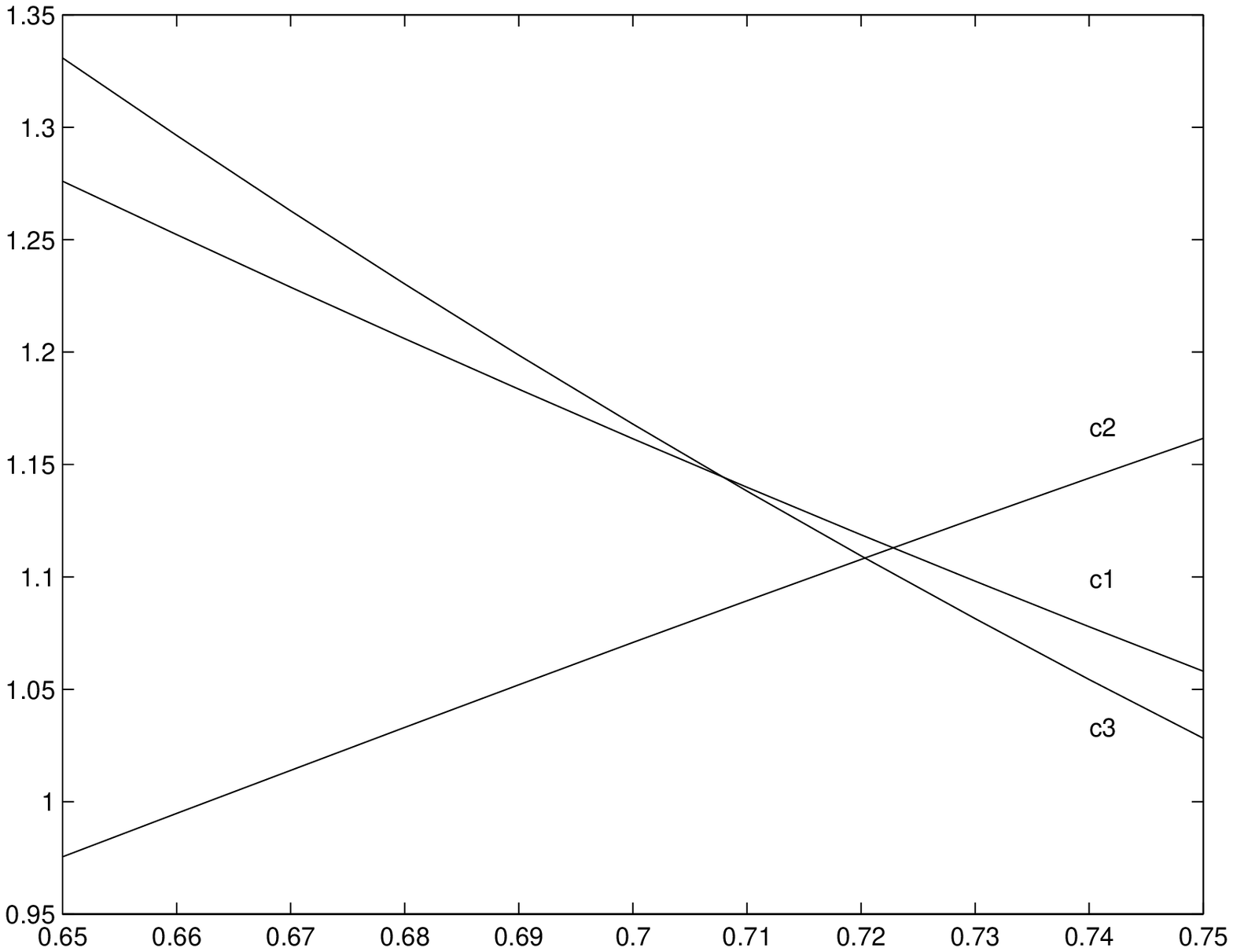}
\epsfxsize 7cm
\epsfbox{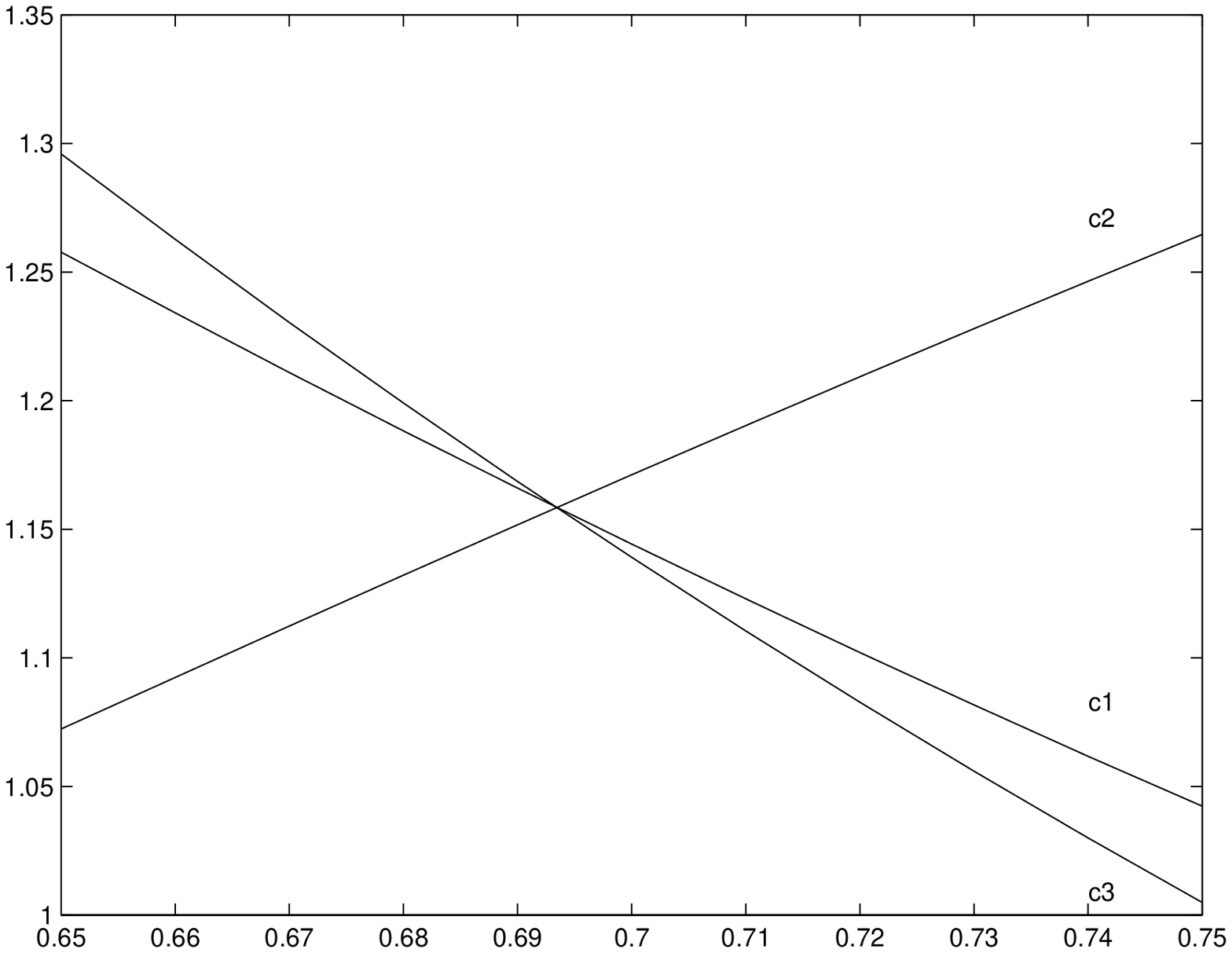}}
\hspace{1.8in}(a)\hspace{2.6in}(b)

\centerline{  
\epsfxsize 7cm   
\epsfbox{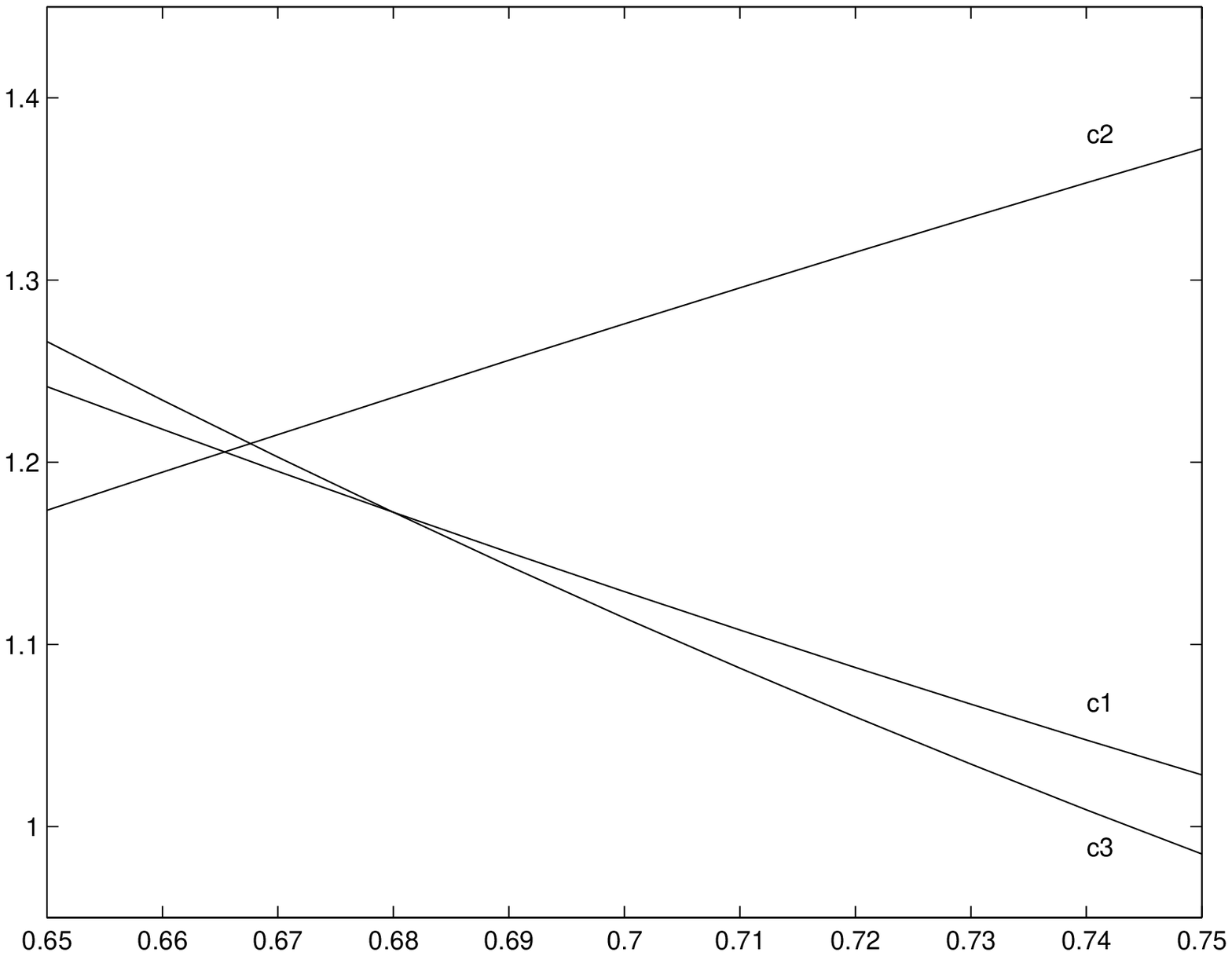}}
\hspace{3.2in}(c) 
\caption{Numeric simulations for (a) $\ld=0.89$ (b) $\ld=0.92$ and (c) $\ld=0.95$.}
\end{figure} 
\eject
\input epsf  
\begin{figure} [ht]  
\centerline{  
\epsfxsize 10cm   
\epsfbox{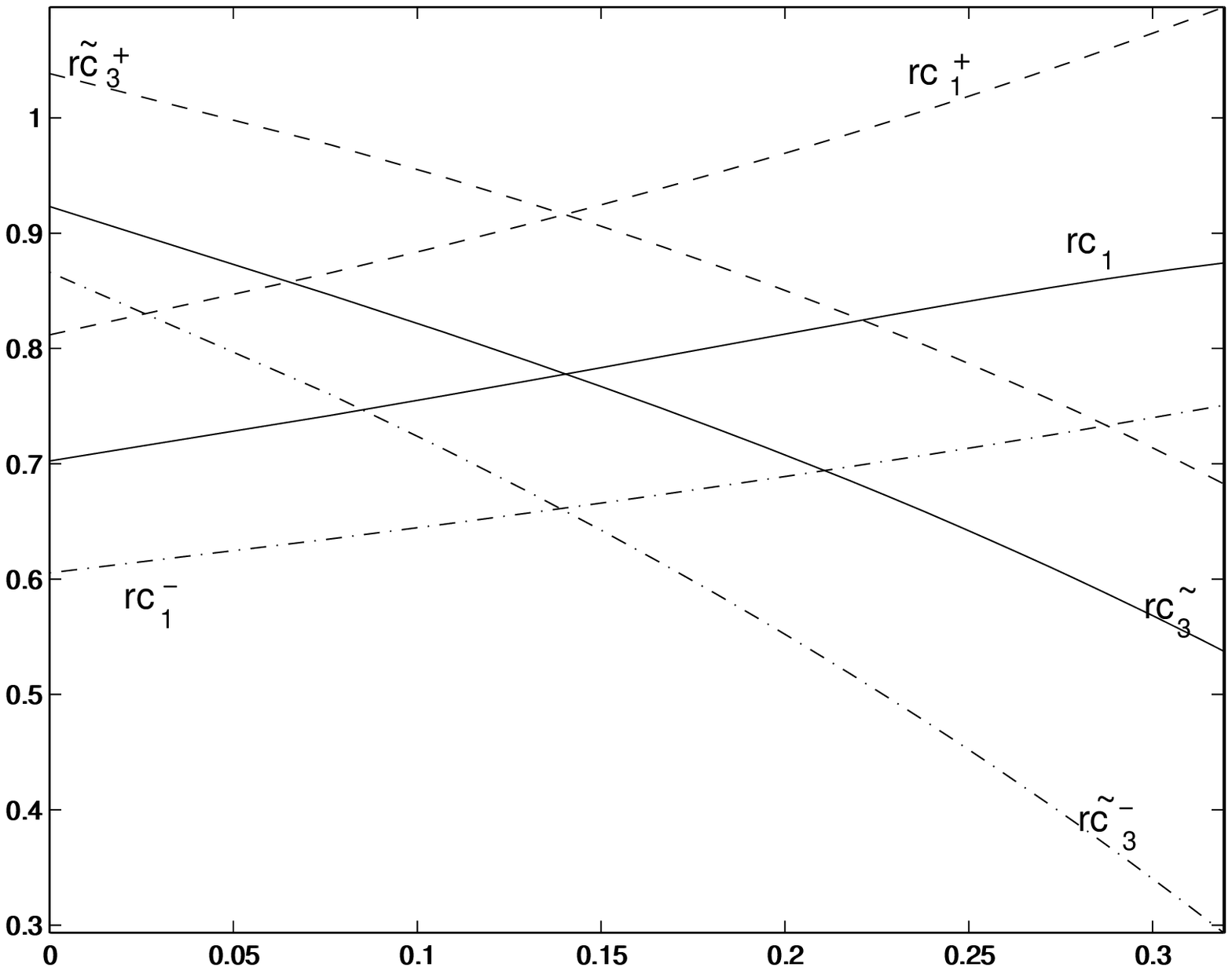}}
\caption{Graphs of $rc_1,rc^{\pm}_{1},r\tc_3$ and $r\tc^{\pm}_{3}$ for $(\rho,\ld)=(-0.05,1.3)$.}
\end{figure}
\ \\
{\bf Proposition 7.1.} \it For any $\rho\in(-0.01,0)$ there exists a point $(\ld^*,\R^*)\in(1,\infty)\times(0,r_0)$ at which the values of $c^2$ in (9) and (15) are positive and coincide simultaneously ($c_1=c_2=c_3$).\rm
\\

The proof of Proposition 7.1 requires some lemmas and estimates we shall deal with in the next section.
\\
\\
{\bf 8. Six main lemmas}
\\
\\
Before starting the demonstration of Proposition 7.1, we need to survey the behaviour of the $x$-roots in (4). This is done in the following lemma:
\\
\\
{\bf Lemma 8.1.} \it For negative $\rho$, $\ld\ge 1$ and $x=a+ib$ one has\rm :\it
\\
(a) The pure imaginary root of (4) is given by $iy$, where $y$ is an increasing function $y=y(r)$. Moreover, $y>r^2+2\srho$ and $y'>2r$ for any positive $r$. 
\\
(b) The function $b(r)$ is negative and decreasing, $b(0)=\srho$ and $b(\infty)=-\ld$.
\\
(c) The function $a(r)$ is positive and $a(0)=\cos\rho$.\rm
\\
\\
\\
{\it Proof}
\\
(a): If $y(r)$ is the pure imaginary root of $(4)$, then 
\BE
   \frac{y(y^2+2y|\srho|+1)}{(y+\ld)^2}=r^2.
\EE

The function $y=y(r)$ is well defined if 
\BE
   y^3+3\ld y^2-(1-4\ld|\srho|)y+\ld>0. 
\EE

The Cardano-Tartaglia discriminant of this last polynomial for $\rho=0$ is $-1/27+2\ld^2/3+\ld^4>0$. Thus (17) holds for $\rho=0$ and consequently for negative $\rho$ as well. Moreover, $y(r)$ is increasing. Now observe the following inequality, which is valid for $r=0$:
\BE
   y>r^2+2\srho.
\EE

If we show that $y'>2r$, then (18) will hold for any positive $r$. From (16)
we have:
\BE
   y'(r)=\frac{1}{r'(y)}=\frac{2r(y+\ld)^3}{y^3+3\ld y^2-(1-4\ld|\srho|)y+\ld}.
\EE 

One easily proves that the right-hand side of (19) is bigger than $2r$. Hence
(18) is valid for every positive $r$.
\\
\\
(b): By setting $x=i\tx$, (4) is equivalent to
\BE
   i(i\tx-iy)(i\tx-a-ib)(i\tx+a-ib)=\tx^3-(2b+y)\tx^2+(a^2+b^2+2by)\tx-(a^2+b^2)y=0.
\EE 

Thus, from (5) and (20) we have $-(2b+y)=R=-2\srho-r^2$ and therefore, $-2b=y-r^2-2\srho$. From item (a) of this Lemma it follows that $b$ and $b'$ are always negative. One easily verifies that $b(0)=\srho$ and $b(\infty)=-\ld$.
\\
\\
(c): The Cardano-Tartaglia discriminant of (4) cannot be negative, for in this case we would get three {\it distinct} pure imaginary roots. This is impossible for $\rho<0$, since (4) is invariant under the map $x\to -\bar{x}$. We still could have $\D=0$ at some $r$, and therefore a double or triple root at this point, all of them pure imaginary. In this case the roots would be $iy$ (simple) and $ib$ (double). But $b$ is negative, according to item (b) of this Lemma, and then (16) could not hold. Therefore, $\D$ is positive. 

Geometrically, this means that the curve $x(r)$ cannot cross the imaginary axis, and consequently $a$ is always positive. It is easy to verify that $a(0)=\cos\rho$ and $a(\infty)=0$. 
\\

\hfill q.e.d.
\ \\

We have just seen that $b=\srho+(r^2-y)/2$. From (5) and (20) we also get $a^2+b^2=\ld^2r^2y^{-1}$. These formulae will be useful in the demonstration of the next lemma, but first define $y_0:=0.2$ and 
\[
   r_0:=\frac{[y_0(y_0^2+2|\srho|y_0+1)]^\m}{y_0+1}.
\]
\\
\\
{\bf Lemma 8.2.} \it For $\ld\ge 1$, $0>\srho\ge -0.01$ and $r\le r_0/\ld$ we have $a_{min}\ge 0.8464$; $a_{max}\le 1.03$; $b_{min}\ge -0.0764$ and $b_{max}=\srho$.\rm 
\\
\\
{\it Proof}. Since $b=\srho+(r^2-y)/2$ and $b'<0$, then it is trivial that $b_{max}=\srho$. Moreover,
\[
   b_{min}\ge -0.01+[r_0^2/\ld^2-y(r_0/\ld)]/2\ge -0.01+[r_0^2-y(r_0)]/2\ge -0.0375\ge -0.0764.
\]

From (5) and (20) we have $a^2=\ld^2r^2y^{-1}-b^2$. Since the derivative of $(y^2+1)/(y+1)^2$ is negative for $y<1$, it follows that 
\[
   a_{min}^2=\frac{y^2+1}{(\frac{y}{\ld}+1)^2}-b_{min}^2
   \ge\frac{y^2+1}{(y+1)^2}-b_{min}^2\ge
   \frac{y_0^2+1}{(y_0+1)^2}-b_{min}^2\ge 0.849^2.
\]

Finally, since $\ld$ can take arbitrarily large values, it is not difficult to see that $a_0^2\le y_0^2+2y_0|\srho|_{max}+1\ge 1.022^2$. 

\hfill q.e.d.
\ \\
\\
\underline{Remark 2}: from now on, we shall apply Lemma 8.2 in several further estimates. Some of them can be slightly improved by making use of $b_{min}=-0.0375$, $a_{min}=0.849$ and $|x|_{min}^2=0.7\dot{2}$.

Now it is easy to evaluate the function $c_1$. We can rewrite (9) as 
\BE
   c^{-1}_1=\frac{2af(\ld)}{(\ld+b)^2+a^2}r-r^2.   
\EE

For $r\in(0,r_0/\ld)$ we apply Lemmas 1 and 2 in order to get the following inequality
\BE
   \frac{2a_{min}f(\ld)}{\ld^2+a_{max}^2}-r<\frac{c^{-1}_1}{r}<
   \frac{2a_{max}f(\ld)}{(\ld+b_{min})^2+a_{min}^2}-r.
\EE

Both the left- and right-hand side of (22) are decreasing with $r$. Moreover, observe (22) and consider $\ld\to\infty$. The bigger $\ld$ is, the ``sooner'' $(rc_1)^{-1}$ gets negative. Because of that, later on we shall work with the variable $\R=\ld r\in(0,r_0)$. 

In order to evaluate the functions $c_2$ and $c_3$, we first recall (15) and define $J_1$ and $J_2$ as the integrals at the numerator and denominator of $c_3$, respectively. In the same way, let $I_1$ and $I_2$ be the numerator and the denominator of $c_2$, respectively. Therefore, we rewrite (15) as 
\BE
   c_2=\frac{I_1}{I_2}\eh\eh{\rm and}\eh\eh
   c_3=\frac{J_1}{\ds\frac{\pi f(\ld)}{(\ld+b)^2+a^2}-J_2}.
\EE   

The next lemma summarises important estimates. It is proved in [RB2] and [RB3]:
\\
\\
{\bf Lemma 8.3.} \it For $t\in(0,1)$ the following equations hold{\rm :} 
\BE
   \int_0^1\frac{2t^2dt}{(t^4+0.1t^2+1)^\m(t^4+0.15t^2+1.06)}+
   \int_0^1\frac{t^\m dt}{(t^2+0.1t+1)^\m(1+0.15t+1.06t^2)}>0.7669;
\EE
\BE
   \int_0^1\frac{2dt}{(t^4+0.1t^2+1)^\m(t^4+0.15t^2+1.06)}+
   \int_0^1\frac{t^{3/2}dt}{(t^2+0.1t+1)^\m(1+0.15t+1.06t^2)}<1.6981;
\EE
\BE
   \int_0^1\frac{2t^2dt}{(t^4+1)^\m(t^4+0.7164)}+
   \int_0^1\frac{t^\m dt}{(t^2+1)^\m(1+0.7164t^2)}>0.9963;
\EE
\BE
   \int_0^1\frac{2dt}{(t^4+1)^\m(t^4+0.7164)}+
   \int_0^1\frac{t^{3/2}dt}{(t^2+1)^\m(1+0.7164t^2)}<2.4499;
\EE
\BE
   \frac{(t^2+1)^\m}{t^2+0.15t+1.06}>-0.3001t+0.94;
\EE
\BE
   \frac{(t^4+1)^\m}{1+0.15t^2+1.06t^4}>-0.37t^2+1;
\EE
\BE
   \frac{(t^2+0.1t+1)^\m}{t^2+0.7164}<-0.69t+1.54;
\EE
\BE
   \frac{(t^4+0.1t^2+1)^\m}{1+0.7164t^4}<1.003;
\EE
\BE
   \int_{0}^{1}\biggl[
               \frac{(t^4-0.1t^2+1)^{-\m}}{(t^2-0.1)^2+0.7164}-
               \frac{t^2(t^4+1)^{-\m}}{1+1.06t^4}
               \biggl](1-t^2)dt<0.764;
\EE
\BE
   \int_0^1\frac{[-0.15-t^2+0.7222(1+t^2+t^4)](t^4-0.1t^2+1)^\m}{[1+1.06t^4][t^4+1.06]}dt>0.443;
\EE
\[
   \int_0^1\biggl[\frac{1}{t^4+1.06}+\frac{t^4}{1+1.06t^4}\biggl]\frac{dt}{(t^4+1)^\m}+
\]
\BE
  -\int_0^1\biggl[\frac{t^2}{(t^2-0.0764)^2+0.7164}+
   \frac{t^2}{(1-0.0764t^2)^2+0.7164t^4}\biggl]\frac{dt}{(t^4-0.1t^2+1)^\m}>0.3294;
\EE
and
\BE
   \sqrt{1-0.1t+t^2}>1-0.0946t+0.473t^2.
\EE\rm
\\

Now we prove the following lemma: 
\\
\\
{\bf Lemma 8.4.} \it Under the same hypothesis of Lemma 8.2, one has 
\BE
   J_{1,min}:=0.7669+1.6981\ld<J_1<J_{1,max}:=0.9963+2.4499\ld\eh\eh{\rm and}
\EE
\BE
   J_{2,min}<J_2<J_{2,max},\eh\eh{\rm where} 
\EE  
\BE
   J_{2,min}=-0.2+(0.6002\ld+1.88)(1-\ld^\m\arctan\ld^{-\m})
             +(0.74\ld^{-1}+2)\ld^{-\m}\arctan\ld^\m-0.74\ld^{-1}
\EE
and
\BE
   J_{2,max}=-0.46+(1.38\ld+3.08)(1-\ld^\m\arctan\ld^{-\m})
             +2.006\ld^{-\m}\arctan\ld^\m.
\EE\rm
{\it Proof}. By splitting the integration interval $(0,\infty)$ into $(0,1]\cup[1,\infty)$ and applying appropriate changes of variable, we finally obtain
\BE
   J_1=\int_0^1\biggl\{\frac{2(t^2+\ld)}{(t^4-2t^2\srho+1)^\m[(t^2-b)^2+a^2]}+
   \frac{t^\m(1+\ld t)}{(t^2-2t\srho+1)^\m[(1-bt)^2+a^2t^2]}\biggl\}dt
\EE
and
\BE
   J_2=\int_0^1\biggl\{\frac{[t(t^2-2t\srho+1)]^\m}{(t+\ld)[(t-b)^2+a^2]}+
   \frac{2(t^4-2t^2\srho+1)^\m}{(1+\ld t^2)[(1-bt^2)^2+a^2t^4]}\biggl\}dt.
\EE 

From Lemma 8.2, it is not difficult to verify the following inequalities:
\[
   J_1>\int_0^1\frac{2(t^2+\ld)dt}{(t^4+0.1t^2+1)^\m(t^4+0.15t^2+1.06)}+
   \int_0^1\frac{t^\m(1+\ld t)dt}{(t^2+0.1t+1)^\m(1+0.15t+1.06t^2)};
\]
\[
   J_1<
   \int_0^1\frac{2(t^2+\ld)dt}{(t^4+1)^\m(t^4+0.7164)}+
   \int_0^1\frac{t^\m(1+\ld t)dt}{(t^2+1)^\m(1+0.7164t^2)};
\]
\[
   J_2>
   \int_0^1\frac{[t(t^2+1)]^\m dt}{(t+\ld)(t^2+0.15t+1.06)}+
   \int_0^1\frac{2(t^4+1)^\m dt}{(1+\ld t^2)(1+0.15t^2+1.06t^4)}
\]
and
\[
   J_2<
   \int_0^1\frac{[t(t^2+0.1t+1)]^\m dt}{(t+\ld)(t^2+0.7164)}+
   \int_0^1\frac{2(t^4+0.1t^2+1)^\m dt}{(1+\ld t^2)(1+0.7164t^4)}.
\]

Now (36) follows from Lemma 8.3. Regarding the inequalities involving $J_2$, 
after we simplify the integrands by Lemma 8.3, we finally obtain (37). 

\hfill q.e.d.
\ \\

At this point we remark that the function $c_3$ has the disadvantage that $a_{min}$, $a_{max}$, $b_{min}$, $b_{max}$, $J_{1,min}$, $J_{1,max}$, $J_{2,min}$ and $J_{2,max}$ will evaluate it by upper and lower functions which do {\it not} depend on $r$. Nevertheless, instead of working with $c_3$ we can work with 
\BE
   \tc_3:=\frac{1}{r}\cdot\biggl(\frac{\pi-2aJ_1r}{2aJ_2-\pi r}\biggl). 
\EE 

While $c_1$ establishes the equality between the residues at the Scherk-ends and $c_3$ the equality between the residue at $z=-i\ld$ and the length of the segment $AB$, the function $\tc_3$ establishes the equality between the residue at $z=x$ and the length of $AB$. In our proof, it will be more convenient to work with $\tc_3$ because $a_{min}$, $a_{max}$, $b_{min}$, $b_{max}$, $J_{1,min}$, $J_{1,max}$, $J_{2,min}$ and $J_{2,max}$ will minimise and maximise (42) by functions which {\it still} depend on $r$. In fact we have
\BE
   \frac{2a_{min}J_{2,min}-\pi r}{\pi-2a_{min}J_{1,min}r}<\frac{\tc_3^{-1}}{r}<
   \frac{2a_{max}J_{2,max}-\pi r}{\pi-2a_{max}J_{1,max}r}.
\EE 

Our problem is to find an $r$ such that $c_1=c_2=\tc_3$. Regarding the equality $c_1=\tc_3$, it is equivalent to
\BE   
   \frac{2aJ_2-\pi r}{\pi-2aJ_1r}=\frac{2af(\ld)}{(\ld+b)^2+a^2}-r.
\EE

By means of Lemmas 8.2 and 8.4, the complexity of (44) was reduced into {\it four} simpler inequalities handling (44) as if the functions there were independent of $r$. In this case, we would have
\BE
    r=\frac{B\pm\sqrt{B^2-4J_1C}}{2J_1},\eh\eh{\rm where}
\EE  
\BE
    B=\frac{2aJ_1f(\ld)}{(\ld+b)^2+a^2}\eh\eh{\rm and}\eh\eh
    C=-J_2+\frac{\pi f(\ld)}{(\ld+b)^2+a^2}.
\EE

If one proves that $C$ is positive and 
\BE
    \frac{B\ld^{-1}}{2J_1}>r_0=\frac{[y_0(y_0^2+2|\srho|y_0+1)]^\m}{y_0+1}\ge r,
\EE
then we have
\BE
    r=\frac{B-\sqrt{B^2-4J_1C}}{2J_1}.
\EE

Hence, (48) establishes the four vertexes of the dashed quadrilateral exemplified in Figure 10. These details will be fulfilled later. Now we are going to analyse the function $c_2$. For convenience of the reader, we rewrite it here:
\BE
    c_2=\frac{\ds\int_{0}^{\infty}\frac{(\ld-t)dt}{f(t)[(t+b)^2+a^2]}}
             {\ds\int_{0}^{\infty}\frac{f(t)dt}{(t-\ld)[(t+b)^2+a^2]}}
       =\frac{I_1}{I_2}.
\EE
\ \\
{\bf Lemma 8.5.} \it Under the same hypothesis of Lemma 8.2, one has 
\BE
   \frac{I_2}{I_1}|_{\ld=1}>0.5798\eh\eh{\rm and}
\EE
\BE
   \limsup_{\ld\to\infty}\ld^{5/2}\frac{I_2}{I_1}|_{\ld=\infty}<0. 
\EE\rm  
{\it Proof}. Once again, by splitting the integration interval $(0,\infty)$ into $(0,1]\cup[1,\infty)$ and applying appropriate changes of variable, we finally obtain 
\BE
   I_1=2\ld^{3/2}\int_{0}^{1}
        \biggl[\frac{(\ld^2t^4+2\ld t^2\srho+1)^{-\m}}
               {(\ld t^2+b)^2+a^2}-
               \frac{t^2(t^4+2\ld t^2\srho+\ld^2)^{-\m}}
               {(\ld+bt^2)^2+a^2t^4}
        \biggl](1-t^2)dt
\EE
and 
\BE
   I_2=2\sqrt{\ld}\int_{0}^{1}
       \biggl[\frac{(t^4+2\ld t^2\srho+\ld^2)^\m}{(\ld+bt^2)^2+a^2t^4}-
              \frac{t^2(\ld^2t^4+2\ld t^2\srho+1)^\m}{(\ld t^2+b)^2+a^2}
       \biggl]\frac{dt}{(1-t^2)}.
\EE

From (52) and Lemma 8.3 one easily computes 
\BE
   I_1|_{\ld=1}<2\int_{0}^{1}
                 \biggl[
                 \frac{(t^4-0.1t^2+1)^{-\m}}{(t^2-0.1)^2+0.7164}-
                 \frac{t^2(t^4+1)^{-\m}}{1+1.06t^4}
                 \biggl](1-t^2)dt<1.528.
\EE

From (53) we have 
\[
   I_2|_{\ld=1}=2\int_{0}^{1}
                 \frac{[(2b-1)t^2+|x|^2(1+t^2+t^4)](t^4+2t^2\srho+1)^\m}
                      {[(1+bt^2)^2+a^2t^4][(t^2+b)^2+a^2]},
\]
hence
\BE
   I_2|_{\ld=1}>2\int_{0}^{1}
                 \frac{[-0.15-t^2+0.7222(1+t^2+t^4)](t^4-0.1t^2+1)^\m}
                      {[1+1.06t^4][t^4+1.06]}.
\EE

By Lemma 8.3, $I_2|_{\ld=1}>0.886$. Thus (50) follows immediately from (54), (55) and the fact that $I_1|_{\ld=1}$ is {\it positive}, which will be proved soon. Regarding (51), from (49) it is immediate to conclude that
\BE
   \liminf_{\ld\to\infty}\ld^{-1}I_1\ge\int_{0}^{\infty}\frac{[t^2+1.06]^{-1}}{f(t)}dt.
\EE
Moreover, from (53) it is not difficult to get 
\BE
   \lim_{\ld\to\infty}\ld^{3/2}I_2=2\srho\int_{0}^{1}(1+t^{-2})dt=-\infty.
\EE

Hence, (51) follows straightforwardly from (56) and (57).
\\

The proof of Lemma 8.5 will be complete if we show that $I_1|_{\ld=1}$ is positive. This follows from 
Lemma 8.6 in the sequence, a stronger result of particular interest. 
\\

\hfill q.e.d.
\ \\
\\
{\bf Lemma 8.6.} \it Under the same hypothesis of Lemma 8.2, $I_1$ remains positive with $r$ in the interval $(0,r_0/\ld)$, for any $(\srho,\ld)\in[-0.01,0)\times[1,\infty)$.\rm
\\
\\
{\it Proof}. First of all, we use (49) to see that
\BE
   I_1=\ld\int_0^\infty\frac{dt}{f(t)[(t+b)^2+a^2]}-\int_0^\infty\frac{tdt}{f(t)[(t+b)^2+a^2]}.
\EE 

By splitting the integration interval $(0,\infty)$ into $(0,1]\cup[1,\infty)$ and applying the changes 
$t\to t^2$ and $t\to t^{-2}$ respectively, we finally obtain
\[
   \frac{I_1}{2}=\ld\int_0^1
                 \biggl[\frac{1}{(t^2+b)^2+a^2}+\frac{t^4}{(1+bt^2)^2+a^2t^4}\biggl]\frac{tdt}{f(t^2)}+
\]
\BE
                   -\int_0^1
                 \biggl[\frac{t^2}{(t^2+b)^2+a^2}+\frac{t^2}{(1+bt^2)^2+a^2t^4}\biggl]\frac{tdt}{f(t^2)}.
\EE

At $\ld=1$ a simple reckon shows that the first integrand of (59) is bigger than the second for any $t\in(0,1)$. Therefore, $I_1|_{\ld=1}$ is positive. Moreover, from Lemma 8.2 and (59) we have
\[
  \frac{I_1}{2}>\ld\int_0^1
                 \biggl[\frac{1}{t^4+1.06}+\frac{t^4}{1+1.06t^4}\biggl]\frac{dt}{(t^4+1)^\m}+
\]
\BE
                   -\int_0^1
                 \biggl[\frac{t^2}{(t^2-0.0764)^2+0.7164}+
                        \frac{t^2}{(1-0.0764t^2)^2+0.7164t^4}\biggl]\frac{dt}{(t^4-0.1t^2+1)^\m}.
\EE 

Because of that, $I_1|_{\ld=1}>0.6588$ by Lemma 8.3. Since the right-hand side of (60) is increasing with $\ld$, it follows the assertion of Lemma 8.6.

\hfill q.e.d.
\ \\

Before we conclude this section, there are two important limits which will be used soon. 
We recall (38) and compute 
\BE
   \lim_{\ld\to\infty}\ld^\m\cdot J_{2,min}=\pi.
\EE

Now observe that  
\[
   1-\ld^\m\arctan\ld^{-\m}=\frac{\ld^{-1}}{3}-\frac{\ld^{-2}}{5}+\frac{\ld^{-3}}{7}-\dots
\]

From (39) it follows that
\BE
   \lim_{\ld\to\infty}\ld^\m\cdot J_{2,max}=1.003\pi.
\EE 

In the next section we finally solve the period problems for the SC surfaces.
\\
\\
{\bf 9. Solution of the period problems}
\\
\\
Let us call $\Q$ the dashed quadrilateral represented in Figure 10. It is determined by the graphs of $rc^{\pm}_1$ and $r\tc^{\pm}_3$, and these functions correspond to the upper and lower estimates in (22) and (43), respectively. The following pictures suggest that $rc_2$ passes through $\Q$ as $\ld$ varies from $1$ to $\infty$.   
\input epsf  
\begin{figure} [ht]  
\centerline{  
\epsfxsize 7.5cm   
\epsfbox{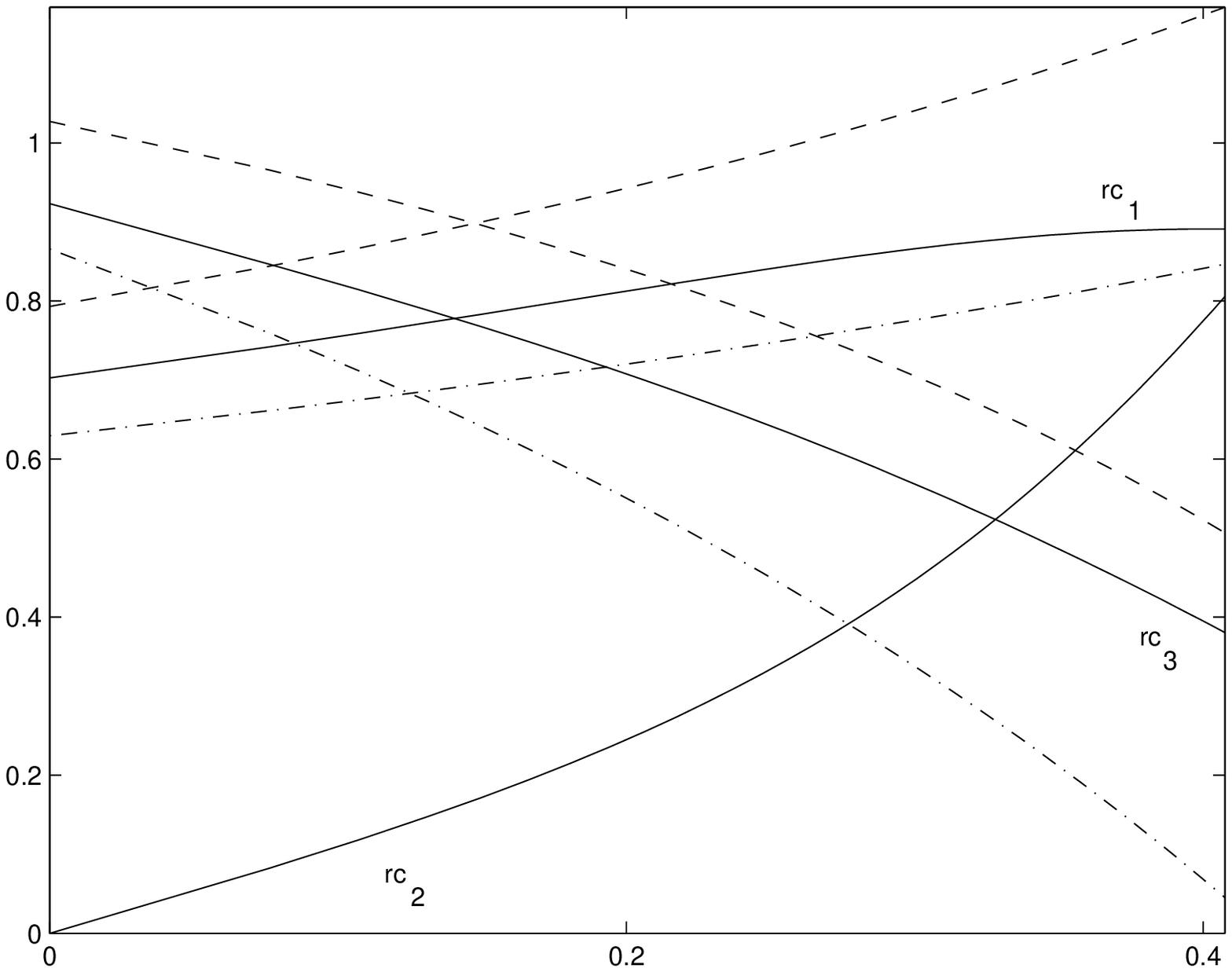}
\epsfxsize 7.2cm
\epsfbox{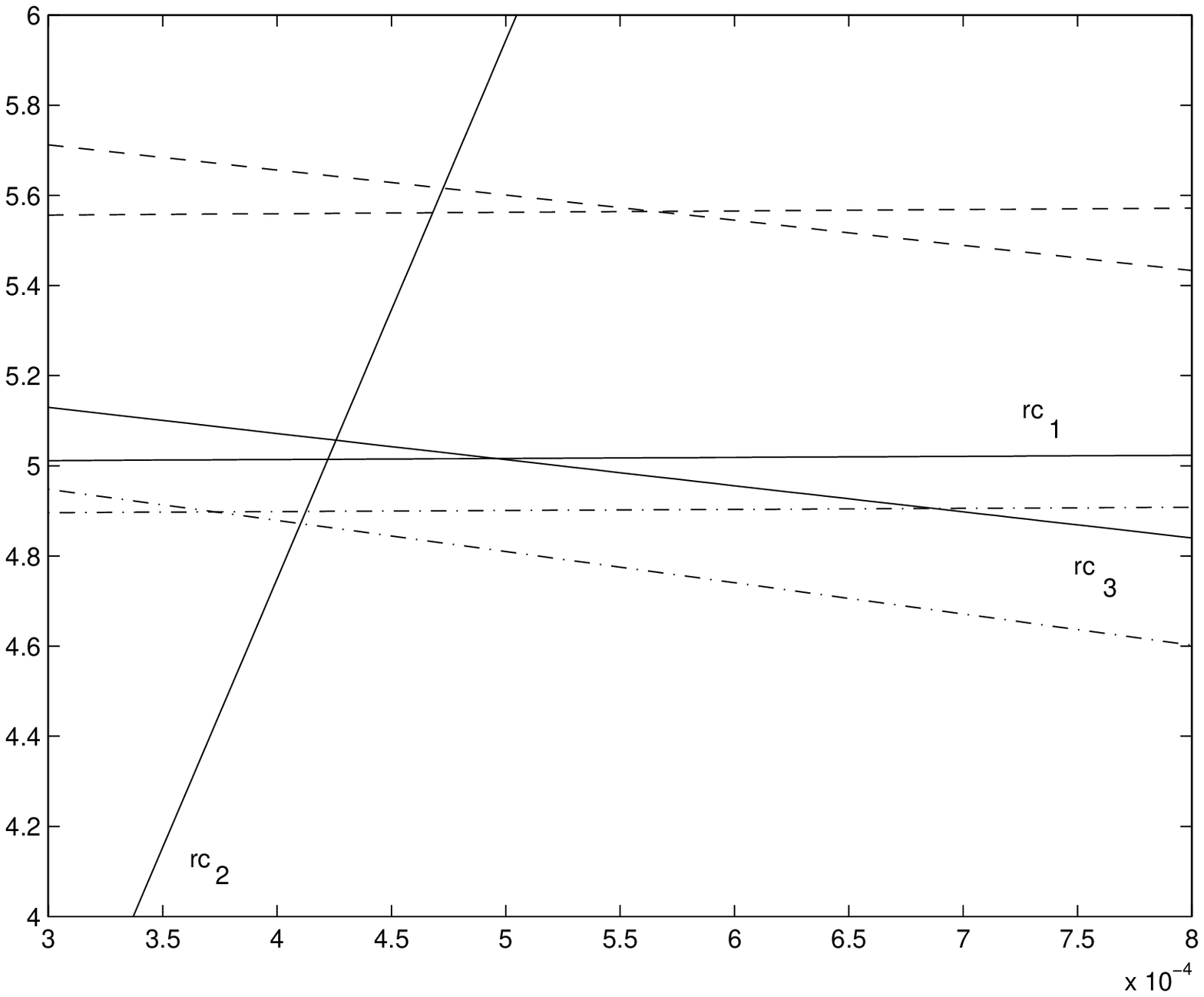}}
\hspace{1.8in}(a)\hspace{2.6in}(b)
\caption{Graphs of $rc_1,rc_2$ and $r\tc_3$ for $\srho=-0.05$ and (a) $\ld=1.3$; (b) $\ld=100$.}
\end{figure}

Our final strategy is to prove the consistency of Figure 11. Without entering into details, one sees that 
Figure 11(a) is true if 
\BE
    rc_2<\max\{rc_1^-,r\tc_3^-\},\eh\eh{\rm for}\eh\eh r\in(0,r_0/\ld),
\EE 
where $r_0\le 0.3808$. From (22) and (23) one sees that (63) will be true if  
\BE
    \frac{2a_{max}f(\ld)}{(\ld+b_{min})^2+a_{min}^2}r-r^2<\frac{I_2}{I_1}.    
\EE

Now we apply Lemmas 8.2 and 8.5 to conclude that (64) holds for $\ld=1$ and any $r$. But since Lemma 8.2 is being used throughout this work, we must always take $r\in(0,r_0/\ld)$. This will complete the first step of our proof. The second step deals with $\ld=\infty$. For technical reasons, we have introduced the variable $\R=\ld r\in(0,r_0)$ in order to not degenerate the interval. By using Lemma 8.2 and Equations (22), (36), (43), (61) and (62), one easily concludes that
\BE
   1.6928<\ld^\m\frac{c_1^{-1}}{r}|_{\ld=\infty}<2.06
\EE
and
\BE
   \frac{1.6928\pi}{\pi-2.874\R}<\ld^\m\frac{\tc_3^{-1}}{r}|_{\ld=\infty}<\frac{2.07\pi}{\pi-5.05\R}.
\EE

From (65) and (66) we see that the quadrilateral $\Q$ will degenerate for $\ld$ big enough, but we still can prove that the graphs of $(rc_1)^{-1}$ and $(r\tc_3)^{-1}$ will always cross. This task is fulfilled later (see Lemma 9.2), and together with $(rc_2)^{-1}|_{\ld=\infty}<0$ (from Lemma 8.5), the period problems will be finally solved. But before going ahead, there are three main questions to be answered:
\\
\\
(a) Does $(r\tc_3)^{-1}$ and $(rc_2)^{-1}$ remain finite for all $(\srho,\ld)\in[-0.01,0)\times[1,\infty)$ 
and $\R\in(0,r_0)$?
\\
(b) Does $(rc_1)^{-1}$ remain positive at all crossings?
\\
(c) Let $(\ld^*(\rho),r^*(\rho))$ be the points at which $c^*(\rho):=c_1=c_2=c_3$. Can we guarantee that 
$0<r^*\sqrt{c^*}\le 1$? 
\\

Condition (c) will avoid intersections at the Scherk-ends. The following lemma will give positive answers 
to all these questions:
\\
\\
{\bf Lemma 9.1.} \it Under the same hypothesis of Lemma 8.2 we have both $(r\tc_3)^{-1}$ and $(rc_2)^{-1}$ 
finite, and $0<r^2c_1<1$ for all $(\R,\srho,\ld)\in(0,r_0)\times[-0.01,0)\times[1,\infty)$.\rm 
\\
\\
{\it Proof}. Consider (22) and define
\BE
   \tr:=\frac{2a_{min}f(\ld)}{\ld^2+a_{max}^2}.         
\EE

It is not difficult to prove that the function $\ld^\m\tr$ is increasing with $\ld$ and $\ld^\m\tr|_{\ld=\infty}=2a_{min}=1.6928$. Moreover, $\tr|_{\ld=1}\ge 1.132>0.3808\ge r_0$. Because of that, $(rc_1)^{-1}$ will remain positive with $r$ in the interval $(0,\ld^{-1}r_0)$, for any $(\srho,\ld)\in[-0.01,0)\times[1,\infty)$.
\\

From (43) we see that $(r\tc_3)^{-1}$ will remain finite whenever
\BE
   r<(2a_{max}J_{1,max})^{-1}\pi.
\EE

By recalling (36), (68) will be true if $\ld r<0.4425$. Hence, the condition $\ld r<r_0$ will satisfy (68) as well. Regarding $(rc_2)^{-1}$, from (49) it will be finite if $I_1\ne 0$, which is exactly the assertion of Lemma 8.6.
\\

It remains to prove that $r^2c_1<1$. From (22) we shall have $r^2c_1<1$ providing
\BE
   r<\frac{2a_{min}f(\ld)}{\ld^2+a_{max}^2}-r,
\EE
which will be satisfied if
\BE
   \ld r<\frac{a_{min}f(1)}{1+a_{max}^2}=0.566\eh\eh\forall\eh(\srho,\ld).
\EE

Therefore, $\ld r<r_0$ will imply (70) and consequently $r^2c_1<1$. 
\\

\hfill q.e.d.
\ \\

Now we prove the following important result:
\\
\\
{\bf Lemma 9.2.} \it At $\R=0$ one has $rc_1<r\tc_3$, and at $\R=r_0$ we have $rc_1>r\tc_3$.\rm 
\\
\\
{\it Proof}. First of all, consider (44) and notice that the inequality $rc_1|_{\R=0}<r\tc_3|_{\R=0}$ is equivalent to 
\BE   
   \frac{2af(\ld)}{\ld^2+a^2}|_{r=0}>\frac{2aJ_2}{\pi}|_{r=0}.
\EE

At $r=0$ we have $(a,b)=(\cos\rho,\srho)$, hence (71) is equivalent to
\BE
   \frac{\pi}{2}\cdot\frac{f(\ld)}{\ld^2+\cos^2\rho}>
   \int_0^1\biggl[\frac{t^2}{t^2+\ld}+\frac{1}{1+\ld t^2}\biggl]\frac{dt}{(t^4-2t^2\srho+1)^\m}.
\EE

By using the inequalities for $t\in(0,1)$
\[
  \frac{t^2}{(1+t^4)^\m}<\frac{\sqrt{2}}{2}t\eh\eh{\rm and}\eh\eh
  \frac{1}{(1+t^4)^\m}<(\frac{\sqrt{2}}{2}-1)t^4+1,
\]
we simplify the integrals in (72) and conclude that it will hold if
\BE
   \frac{\pi\sqrt{1-0.1\ld^{-1}+\ld^{-2}}}{2(1+\ld^{-2})}>
   \frac{\sqrt{2}}{4}\ld^\m\ln(1+\ld^{-1})+[(\frac{\sqrt{2}}{2}-1)\ld^{-2}+1]\arctan\ld^\m+
   (1-\frac{\sqrt{2}}{2})\ld^{-\m}(\ld^{-1}-\frac{1}{3}).
\EE

By setting $\Ld:=\ld^{-1/2}$ and using Lemma 8.3, one sees that (73) will be true if
\BE
   \frac{\pi(1-0.0946\Ld^2+0.473\Ld^4)}{2(1+\Ld^4)}>
   \frac{\sqrt{2}}{4}\frac{\ln(1+\Ld^2)}{\Ld}+[(\frac{\sqrt{2}}{2}-1)\Ld^4+1]\arctan\frac{1}{\Ld}+
   (\frac{\sqrt{2}}{2}-1)\Ld(\frac{1}{3}-\Ld^2).
\EE

Since $\ln(1+\Ld^2)<\Ld^2-\Ld^4/4+\Ld^6/3$ and $4\arctan\Ld^{-1}<\pi(2-\Ld)$, (74) will hold providing 
\[
   [3\pi(\sqrt{2}-2)-2\sqrt{2}]\frac{\Ld^8}{24}+\frac{\pi}{4}(2-\sqrt{2})\Ld^7+(\frac{5\sqrt{2}}{8}-1)\Ld^6
   +(\frac{\pi\sqrt{2}}{8}-\frac{\sqrt{2}}{2}+\frac{1}{3})\Ld^4+
\]
\BE
    +\frac{\pi}{4}(0.946-\sqrt{2})\Ld^3+(\frac{5\sqrt{2}}{8}-1)\Ld^2-0.0473\pi\Ld+
    \frac{1}{3}+\frac{\pi}{4}-\frac{5\sqrt{2}}{12}>0.
\EE

Now (75) will be valid if
\BE
    F_1:=-0.348\Ld^8+0.46\Ld^7-0.12\Ld^6+0.181\Ld^4-0.368\Ld^3-0.12\Ld^2-0.15\Ld+0.465>0,
\EE 
and since $F_2:=[F_1/(1-\Ld)-0.465]/\Ld=0.348\Ld^6-0.112\Ld^5+0.008\Ld^4+0.008\Ld^3-0.173\Ld^3+0.195\Ld^2+0.315$, we can simply verify whether $F_2$ is positive in $[0,1]$. But this comes directly from $0.112<0.195$ and $0.173<0.315$. 

We have just established the first inequality of Lemma 9.2. Now we deal with the second and fulfil the details briefly explained in (44)-(48). From (21) and (42), we recall (44) and observe that $rc_1|_{\R=r_0}>r\tc_3|_{\R=r_0}$ if
\BE   
   \frac{2aJ_2-\pi r}{\pi-2aJ_1r}|_{r=\frac{r_0}{\ld}}>\frac{2af(\ld)}{(\ld+b)^2+a^2}-r|_{r=\frac{r_0}{\ld}}.
\EE

Condition (77) only makes sense for a non-zero denominator at its left-hand side. But $\pi-2aJ_1r>\pi-2a_{max}J_{1,max}r$, and (36) implies 
\BE
   \frac{\pi}{2a_{max}(0.9963\ld^{-1}+2.4499)}\ge 0.4425>r_0.
\EE
Hence, (78) allows us to work with (77). Now fix $r=r_0/\ld$ and observe that (77) will hold if
\[   
   \frac{2aJ_2-\pi r}{\pi-2aJ_{1,min}r}>\frac{2af(\ld)}{(\ld+b)^2+a^2}-r,
\]
or equivalently
\[
   (J_2-J_{1,min}r^2)[(\ld+b)^2+a^2]>[\pi-2aJ_{1,min}r]f(\ld),
\]
which in its turn will be valid providing
\BE
   (J_{2,min}-J_{1,min}r^2)[(\ld+b_{min})^2+a_{min}^2]>[\pi-2a_{min}J_{1,min}r]\sqrt{\ld(\ld^2+1)}.
\EE

According to the remark at Lemma 8.2, we might have been using $|x|_{min}^2=0.7\dot{2}$, $a_{min}=0.849$ and $b_{min}=-0.0375$. Moreover, $0.38<r_0<0.3808$. Therefore, at $\ld=1$ the left- and right-hand sides of (79) take the approximate values $2.19$ and $2.28$, respectively. Hence (79) holds at $\ld=1$. We recall $\Ld=\ld^{-1/2}$ and observe the following inequalities, valid for $\Ld\in(0,1)$:
\[
   4\arctan\Ld^{-1}>\Ld^2-4\Ld+3+\pi,
\]
\[
   \arctan\Ld<(\pi/4-2/3)\Ld^4-\Ld^3/3+\Ld\eh\eh{\rm and}\eh\eh
\]
\BE
   32\sqrt{1+\Ld^4}<\sqrt{2}(23+10\Ld^4-\Ld^8).
\EE

One easily verifies (80) by taking derivatives and handling polynomials with roots in $\real\setminus(0,1)$. Consider (36), (38), (79) and (80). Since $0.38^2<r_0^2\le 0.145$, one sees that (79) will hold providing
\[
   \{-0.145\Ld^2(0.7669\Ld^2+1.6981)-0.2-0.74\Ld^2+\frac{\Ld}{4}(\Ld^2-4\Ld+3+\pi)(2+0.74\Ld^2)
\]
\[
   -(0.6002+1.88\Ld^2)[(\pi/4-2/3)\Ld-1/3]\}(\Ld^{-4}-0.075\Ld^{-2}+0.7\dot{2})
\]
\BE
   >\frac{\sqrt{2}}{32}\Ld^{-3}[\pi-0.64524(0.7669\Ld^2+1.6981)](23+10\Ld^4-\Ld^8).
\EE  

After rearranging terms in (81) we see that it will hold if
\[
   F_3:=-0.022\Ld^{10}+0.22\Ld^8-0.615\Ld^7+1.225\Ld^6-1.641\Ld^5
\]
\[
   +1.34\Ld^4-0.675\Ld^3+1.69\Ld^2-2.36\Ld+0.9199>0,
\]
which in its turn will be valid providing $F_4:=F_3-0.022(\Ld^8-\Ld^{10})-0.0039>0$. But $F_4=F_5(\Ld-0.786)^2+F_6$, with $F_6>0$ in $(0,1)$ and $F_5>F_7:=0.198\Ld^6-0.304\Ld^5+0.6\Ld^4-0.5\Ld^3+0.2\Ld^2-0.05\Ld+1.48$. Since $1.48>0.05+0.5+0.304$, then $F_7$ is positive in $[0,1]$. Therefore, (77) holds for any $\ld\in[1,\infty)$ and consequently $rc_1>r\tc_3$ at $\R=r_0$. 
\\

\hfill q.e.d.
\ \\

Now we are ready to prove Proposition 7.1. For convenience of the reader, we restate it here as
\\
\\  
{\bf Proposition 9.1.} \it For any $\rho\in(-0.01,0)$ there exists a point $(\ld^*,r^*)\in(1,\infty)\times(0,r_0)$ at which the values of $c^2$ in (9) and (15) are positive and coincide simultaneously ($c_1=c_2=c_3$).\rm
\\
\\
{\it Proof}. Lemmas 9.1 and 9.2 guarantee the existence of a curve $\be:[0,1)\to[1,\infty)\times[0,r_0]$, given by $\be(s)=(\ld(s),\R(s))$, such that $c_1(\rho,\ld(s),\R(s))=\tc_3(\rho,\ld(s),\R(s))$ for any $s\in[0,1)$, with fixed $\srho\in[-0.01,0)$. Moreover, if $\CC:=\be([0,1))\subset[1,\infty)\times[0,r_0]$ then $c_1|_{\CC}>0$ by Lemma 9.1. 
\\

As we briefly discussed at the beginning of this section, Lemmas 8.2 and 8.5 assure that the inequality
\BE
    \frac{2a_{max}f(\ld)}{(\ld+b_{min})^2+a_{min}^2}r-r^2<\frac{I_2}{I_1}   
\EE 
holds for $\ld=1$ and $r\in(0,r_0)$. Since (82) implies 
\[
    (rc_2)^{-1}>(rc_1^-)^{-1}\ge(rc_1)^{-1}\eh\eh{\rm for}\eh\eh\ld=1\eh\eh{\rm and}\eh\eh r\in(0,r_0),
\] 
we have that 
\BE
    (rc_2)^{-1}|_{\be(0)}>(rc_1)^{-1}|_{\be(0)}.
\EE

Equation (83) is valid because $\be(0)=(\ld(0),\R(0))=(1,\R(0))$. We cannot guarantee that $\Lim{s\to 1}{\ld(s)}=\infty$, but of course there is a sequence $\{s_n\}_{n\in\natl}\subset[0,1)$ with $\Lim{n\to\infty}{s_n}=1$ and $\Lim{n\to\infty}{\ld(s_n)}=\infty$ (see Figure 12). 
\\
\input epsf 
\begin{figure} [ht] 
\centerline{ 
\epsfxsize 7cm  
\epsfbox{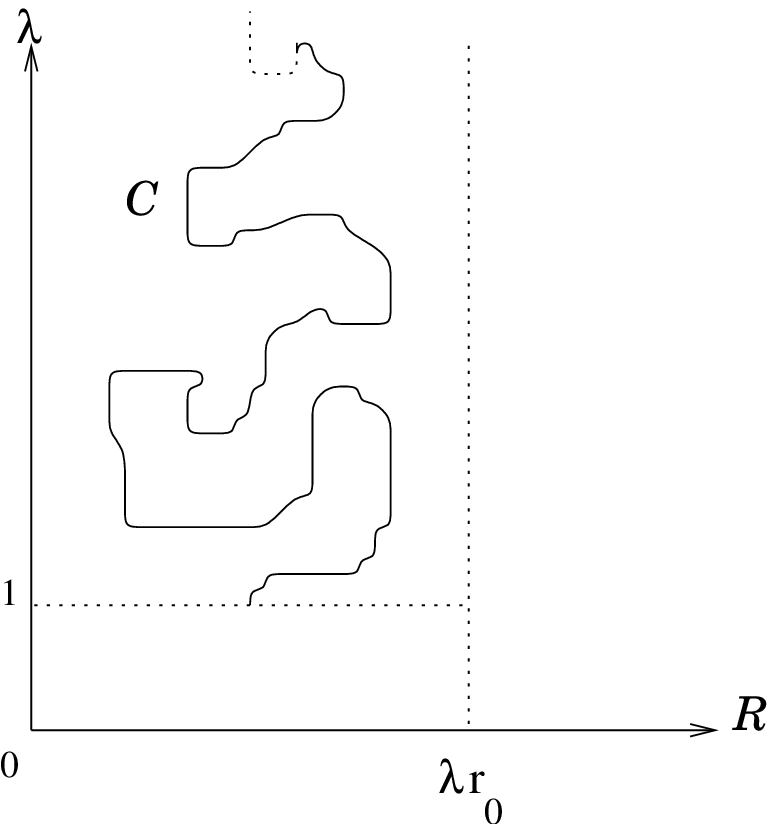}}  
\caption{A scheme of the curve $\CC=\be\big([0,1)\big)\in(0,r_0)\times[1,\infty)$.}   
\end{figure} 

From Lemma 8.5, there exists $m\in\natl$ such that 
\BE
    (rc_2)^{-1}|_{\be(s_m)}\le 0<(rc_1)^{-1}|_{\be(s_m)}.
\EE

Now, (83) and (84) imply the existence of an $s^*\in[0,1)$ such that $c_2|_{\be(s^*)}=c_1|_{\be(s^*)}$. Since $c_1|_{\be(s^*)}=\tc_3|_{\be(s^*)}$ and $c_1$ is positive, there exists $(\ld^*,r^*)=(\ld(s^*),\R(s^*)/\ld(s^*))$ such that $c_1=c_2=\tc_3$. The equality $c_1=\tc_3$ implies $c_1=c_3$. We can define $c^*(\rho):=c_1=c_2=c_3$, where $r^*\sqrt{c^*(\rho)}\in(0,1)$ by Lemma 9.1. This value depends only on $\rho\in(-0.01,0)$.

\hfill q.e.d.    
\\

Now we have demonstrated items (a), (b) and (c) from Theorem 1.1. The next and last section deals with item (d), which will finally conclude our work. 
\\
\\
{\bf 10. The embeddedness of the Scherk-Costa surfaces} 
\\
\\
This final section is devoted to the proof that the SC surfaces are embedded. In order to do this, we shall use a technique similar to that one develop in [FM], where the authors obtained a new embeddedness proof for the Hoffman-Karcher-Wei singly periodic genus-one helicoid (see [HKW]). We shall construct a continuous one-parameter family of minimal surfaces with boundary, starting at an ``open half-Costa surface" (see Figure 2(b)) and ending at a fundamental piece of our surface contained in the upper half-space. We proceed as follows.

For each $s\in[0,1)$ define $r(s):=r^*(1-s)$, $\rho(s):=\rho^* (1-s)$ and $\ld(s):=\ld^*(1-s)+s[2r(s)\mu]^{-\frac23}$, where $r^*$, $\rho^*$ and $\ld^*$ are the parameters that close the period problem for one of our surfaces. We take $x(s)$ given by (4) and 
$$
   c(s):=\left\{\begin{array}{l}\min\{c_1(s),c_2(s)\}\quad\mbox{if}\eh c_2(s)>0,\\ \\
                               c_1(s)\quad\mbox{if}\eh c_2(s)<0, 
                \end{array}\right.
$$
where $c_1(s)$ and $c_2(s)$ are defined by (9) and the first equation of (15), respectively. 
 
Let $T_s$ be the torus with algebraic equation $z'^2=-iz(z-e^{i\rho(s)})(z+e^{-i\rho(s)})$. Consider the twice punctured disk in $T_s$, 
$$ 
   M_s=z^{-1}(\{\zeta\in\C\::\:Re(\zeta)>0\}\setminus\{x\}),
$$
represented in Figure 13. On $M_s$ we have the Weierstrass data given by (2) and (3): 
\begin{eqnarray}
   g_s  &=& c(s)\frac{\sqrt{-iz(z-e^{i\rho(s)})(z+e^{-i\rho(s)})}}{z+\ld(s) i},\label{eq:gs}\\
   dh_s &=& \frac{dz}{(z-x(s))(z+\overline{x(s)})}.\label{eq:dhs}
\end{eqnarray}

These meromorphic data define a conformal, minimal, multi-valued immersion $X_s$ of $M_s$ into $\real^3$. Namely, they define a conformal minimal immersion $\widehat X_s: \widehat M_s \rightarrow \real^3$, where $\widehat M_s$ is the universal covering of $M_s$. The surface $\widehat X_s(\widehat M_s)$ is invariant under a horizontal translation corresponding to the period $\vec p_{x(s)}$ associated to the end $x(s)$. 

As a matter of fact, there are intermediate coverings $\widehat{\pi}:\widehat M_s\to\widetilde M_s$ and $\widetilde{\pi}:\widetilde M_s\to M_s$ such that $\widetilde{\pi}\circ\widehat{\pi}$ is universal, and $X_s\circ\widetilde{\pi}$ is {\it already} univalent. In order to describe $\widetilde{\pi}$, take a loop $\LL$ in $M_s$ which bounds a disk containing $x(s)_+$ and $x(s)_-$. Now $\langle[\LL]\rangle$ is an infinite cyclic subgroup of $\pi_1(M_s)$ and from [B,\S III.8] we have a conformal covering $\widetilde{\pi}:\widetilde M_s\to M_s$, where $\widetilde M_s$ is a Riemann surface with a freely acting infinite cyclic group $\mathbf{P}$ of automorphisms. The group $\mathbf{P}$ is completely determined by $\langle[\LL]\rangle$. Moreover, $M_s=\widetilde M_s/\mathbf{P}$.

Let $\widetilde{X}_s:\widetilde M_s\longrightarrow\real^3$ be the univalent conformal minimal immersion.  

From (7) we know that $\vec p_{x(s)}=(0,k_2(s),0).$ Taking the definition of $c(s)$ into account one has:
\begin{enumerate}[(a)] 
\item there are no contact between the straight half-lines of $\deh\widetilde X_s(\widetilde M_s)$, except for $A=A'$ when $c(s)=c_2(s)$, as illustrated in Figure 14;
\item the period $\vec p_{x(s)}$ is longer than (or equal to) the period $\vec{p}_{\ld(s)}$ associated to $-i\ld$ (see Figure 14). 
\end{enumerate}

Our first step will consist of studying the behaviour of the Gauss map as $s \to 1$. This is the purpose of the following lemma:
\\
\\
{\bf Lemma 10.1.} \it There exist real numbers $\kappa>0$ and $s_0\in(0,1)$ such that, for any $s\in(s_0,1)$
\begin{enumerate}[(i)]
\item $|g_s|<1$ in $z^{-1}(D(1,1/\kappa))$ and
\item $|g_s|>1$ in $z^{-1}(\{\zeta\in\C\;:\;Re(\zeta)\geq 0,\;Im(\zeta)<-\kappa\})$.
\end{enumerate}\rm
\ \\
{\it Proof}. Note that $\ds\lim_{s\to 1}x(s)=1$. Hence, it is not hard to see that $\ds\lim_{s\to 1}\frac{c^2(s)}{\ld^2(s)}=2\mu^2$. Moreover, it is clear that $g_s^2$ can be viewed as a function on the $z$-plane. We then have a family of meromorphic functions $\{g_s^2\}_{s\in[0,1)}$ on $\C$ which converges for $s\to 1$ uniformly on compact sets of $\C$ to $g_0^2=2i\mu^2 z(z^2-1)$. Hence, if $\kappa$ is big enough and $s_0$ is sufficiently close to $1$, then {\it (i)} and {\it (ii)} will hold.

\hfill q.e.d.
\ \\

Take a set of Cartesian co-ordinates in $\real^3$  such that the symmetry of $\widetilde X_s(\widetilde M_s)$ induced by the elliptic involution $(z,w)\mapsto(z,-w)$ coincides with the $180^\circ$ rotation about the $x_3$-axis. Given $\Rg>0$, we define in $\real^3$ the set ${\cal B}_\Rg:=\bigcup_{n\in\Z}B(n\cdot\vec p(x(s)),\Rg)$.
\\
\\
{\bf Lemma 10.2.} \it There exist $R_0>0$ and $s_1\in(0,1)$ such that $\widetilde X_s(\widetilde M_s)\setminus\overline{\cal B}_{R_0}$ is embedded, $\forall s\in(s_1,1)$.\rm
\\
\\
{\it Proof}. As our family of minimal disks converges in compact sets of $\real^3$ to an open Costa surface, then we can assume that $\ds\widetilde{\pi}\left(\widetilde X_s^{-1}\left(\widetilde X_s( \widetilde M_s)\setminus\overline{\cal B}_\Rg\right)\right)$ has exactly three connected components $\Omega_1$, $\Omega_2,$ and $\Omega_3$. The set $\Omega_1$ contains the ends $x_+$ and $x_-$, $\Omega_2$ contains the horizontal end $i\ld_+$ and $\Omega_3$ contains $i\ld_-$. We write $\widetilde{\Omega}_k=\widetilde{\pi}^{-1}(\Omega_k)$, $k=1,2,3$. Our purpose is to prove that the $\widetilde X_s(\widetilde\Omega_k)$ are disjoint graphs on the $(x_1,x_2)$-plane.
\input epsf 
\begin{figure} [ht] 
\centerline{ 
\epsfxsize 9cm  
\epsfbox{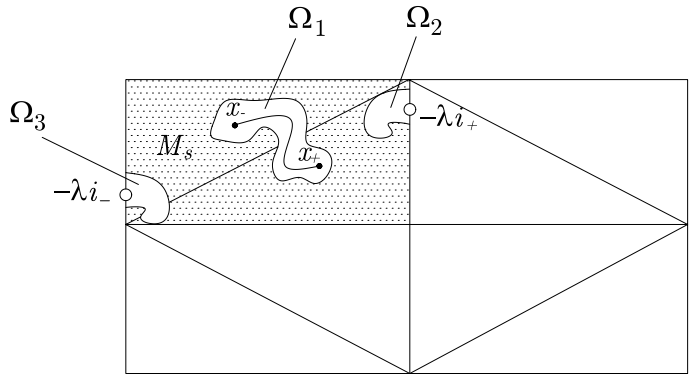}}
\caption{The sets $\Omega_i$, $i=1,2,3$.}
\end{figure}

We first consider $\Omega_1$.  Let $\kappa>0$ as in Lemma 10.1 and take $s_1>s_0$. If $s_1$ is close enough to $1$, we can assume that $\Omega_1\subset z^{-1}(D(1,1/\kappa))$. Using once again the fact that our family converges in each $\overline{\cal B}_\Rg$ to an open Costa surface, then we guarantee (for $s_1$ close to $1$ and $\Rg$ sufficiently large) that $\widetilde X_s(\deh\widetilde\Omega_1)$ is a set of pairwise disjoint Jordan curves. Only one of these curves intersects the plane $\{x_2=0\}$, and in exactly two points: $P_-=(-k_1,0,k_3)$ and $P_+=(k_1,0,k_3)$, $k_1,k_3>0$. We want to study the set $\G=\widetilde X_s(\widetilde \Omega_1) \cap \{x_2=0\}$. It is known from [Ch] that this intersection is a set of analytic curves. Regarding $\G$, we can say that:
\begin{enumerate}[(a)]
\item Only one curve of $\G$ arrives at $P_-$ and only another one at $P_+$. Furthermore the slope of these curves is negative at $P_-$ and positive at $P_+$. This is also a consequence of the convergence of our surfaces to the open Costa surface.
\item The curves in $\G$ must diverge to the upper Scherk ends $x_+$ and $x_-$, since the immersion is proper. Moreover, there is only one curve in $\G$ diverging to each end (the ends are embedded).
\item All curves in $\G$ are graphs over the $x_1$-axis, because of Lemma 10.1 and the fact that $\Omega_1\subset z^{-1}(D(1,1/\kappa))$. In particular, there are no compact curves in $\G$.
\end{enumerate} 
These properties imply that $\G$ consists of two curves, $\G_-$ and $\G_+$. The curve $\G_-$ starts at $P_-$ and diverges in the negative $x_1$-axis direction. It happens the analogous to $\G_+$. In particular, we have proved that the third projection $\pg_3|_{\widetilde \Omega_1}$ is not onto. From our assumptions, we know that $\pg_3(\deh\widetilde{\Omega}_1)$ is a sequence of Jordan curves $\{\gamma_j\}_{j\in\Z}$ in the $(x_1,x_2)$-plane satisfying $\gamma_{j+1}=\gamma_j+\vec p_{x(s)}$, for any $j\in\Z$. Now define 
$$
   U_1:=\bigcup_{j\in\Z}\mbox{Int}(\gamma_j),\qquad U_2:=\bigcap_{j\in\Z}\mbox{Ext}(\gamma_j), 
$$ 
where $\mbox{Int}(\gamma_j)$ and  $\mbox{Ext}(\gamma_j)$ are the interior and exterior of $\gamma_j$, respectively. 
\input epsf 
\begin{figure} [ht] 
\centerline{ 
\epsfxsize 9cm  
\epsfbox{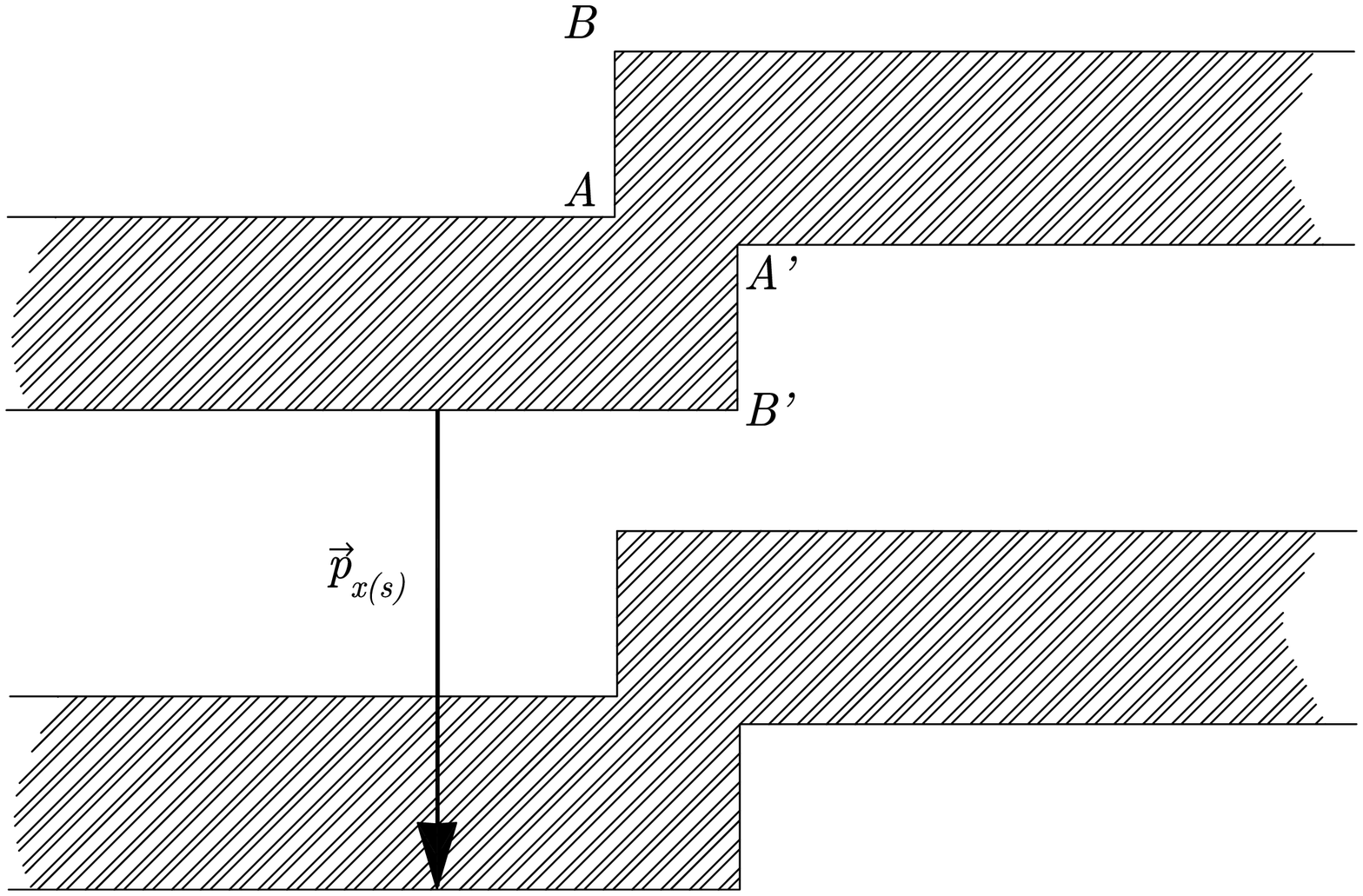}}
\caption{The boundary of $\widetilde{X}_s(\widetilde{M}_s).$}
\end{figure} 

Since $|g_s|<1$ in $\widetilde{\Omega}_1$, then ${\pg_3}|_{\widetilde{\Omega}_1}:\widetilde{\Omega}_1\rightarrow \{x_3=0\}$ is a local diffeomorphism. Hence $\pg_3(\widetilde{\Omega}_1)$ is an open subset of the plane $\{x_3=0\}$. Moreover, ${\pg_3}|_{\widetilde{\Omega}_1}$ is proper. Then $\pg_3(\widetilde\Omega_1)$ is closed in $\{x_3=0\}\setminus\pg_3(\deh\widetilde{\Omega}_1)$. Therefore, either $U_1 \subset\pg_3(\widetilde{\Omega}_1)$ or $U_1\cap\pg_3(\widetilde{\Omega}_1)=\emptyset$. But ${\pg_3}|_{\widetilde{\Omega}_1}$ is not onto, so $\pg_3(\widetilde{\Omega}_1)=U_2$. This implies that ${\pg_3}|_{\widetilde{\Omega}_1}:\widetilde{\Omega}_1 \longrightarrow U_2$ is a covering map. Since ${\pg_3}|_{\widetilde{\Omega}_1}$ is one-to-one at the ends, then ${\pg_3}|_{\widetilde{\Omega}_1}$ is a diffeomorphism, namely $\widetilde{\Omega}_1$ is a graph.

In order to prove that $\widetilde\Omega_2$ and $\widetilde\Omega_3$ are graphs one proceeds in a similar way as for $\widetilde \Omega_1$, but now uses {\it (ii)} in Lemma 10.1 and applies similar arguments about covering maps.

\hfill q.e.d.
\\
\\
{\bf Lemma 10.3.} \it There exists $s_2\in(0,1)$ such that $\widetilde X_s: \widetilde M_s \longrightarrow \real^3$ is an embedding, for all $s\in(s_2,1)$.\rm 
\\
\\
{\it Proof}. By using Lemma 10.2 we have obtained a radius $R_0>0$ and an interval $(s_1,1)$ such that $\widetilde X_s(\widetilde M_s)\setminus\overline{\cal B}_{R_0}$ is embedded. In the balls $\overline{B}(n\cdot\vec p(x(s)),R_0)$, $n\in\Z$, $\widetilde X_s$ converges uniformly to an open half-Costa surface. Thus we can use the {\it interior maximum principle} from [Sn] and get $s_2\in(0,1)$ with $\widetilde X_s$ an embedding for all $s\in(s_2,1)$.

\hfill q.e.d.

At this point we can prove the main result of this last section.
\\
\\
{\bf Theorem 10.1.} \it Consider $\rho\in(-0.01,0)$ and let $\ld^*$, $r^*$ be the parameters given by Proposition 9.1. Then the corresponding SC surface is embedded.\rm
\\
\\
{\it Proof}. As in the proof of Lemma 10.2, we denote $\widetilde{X}_s:\widetilde{M}_s\longrightarrow\real^3$ the minimal immersion (with boundary) obtained by successive translations of $M_s$ in the period $\vec p_{x(s)}$. Define the set 
$$
   {\cal I}=\{s\in[0,1)\::\:\widetilde{X}_s:\widetilde{M}_s\longrightarrow\real^3\quad
   \mbox{is an embedding}\}.
$$

From Lemma 10.3 one has $(s_2,1)\subset{\cal I}$, in particular ${\cal I} \neq \emptyset$.
 
Consider $s'\in {\cal I}$. Taking into account our construction of the Weierstrass representation of $X_s$, one has:
\begin{enumerate}[(a)]
\item The ends $x(s)_+$, $x(s)_-$, $-i \lambda(s)_+$ and $-i \lambda(s)_-$ do not intersect because $0<|g_s(x(s)_\pm)|<1$, for $s\in[0,1)$. This comes from the inequality $0<c(s)\leq c_1(s)$ and Lemma 9.1, from which one has $0<r^2c_1<1$. Now use the fact that $g_s(x(s))=\pm r(s)\sqrt{c(s)}$. 
\item The period $\vec p_{x(s)}$ is longer than $\vec p_{\ld(s)}$. This is due to the choice we made of the L\'opez-Ros parameter $c(s)$.
\item The boundary of $\widetilde X_s(\widetilde M_s)$ is a sequence of parallel non-compact polygonal lines. Each connected component consists of two half-lines parallel to the $x_1$-axis joined by a segment which is parallel to the $x_2$-axis. The boundary is described in Figure 14. The choice of $c(s)$ also guarantees that the distance between the connected components is positive $\forall s \in \: [0,1)$. 
\end{enumerate}
Then, there exist $\epsilon,\eta>0$ such that $\widetilde{X}_s(\widetilde{M}_s)\cap(\real^3\setminus{\cal C}_\eta)$ is embedded, $\forall s\in(s'-\epsilon,s'+\epsilon)$, where ${\cal C}_\eta=\{(x_1,x_2,x_3)\in\real^3\::\:x_1^2+x_3^2<\eta^2\}$. If $\widetilde{X}_s$ were not injective for some $s\in(s'-\epsilon,s'+\epsilon)$, then self-intersections of $\widetilde{X}_s(\widetilde{M}_s)$ would be in ${\cal C}_\eta$. As the length of the period $\vec p_{x(s)}$ is bounded, then the first contact must occur in a compact region of ${\cal C}_\eta$. Moreover, as $\widetilde X_s(\widetilde M_s)$ is contained in the upper half-space, there are no contacts of interior points with points of the boundary, this latter contained in the plane $\{x_3=0\}$. So, there would be contact between interior points, contradicting the {\it classical maximum principle} (see [Sn]). Hence $(s'-\epsilon,s'+\epsilon)\subset{\cal I}$, which implies that $\cal I$ is open.

Now take a sequence $\{s_n\}_{n\in\natl}$ in $\cal I$ converging to $s''\in[0,1)$. Assume that $\widetilde{X}_{s''}$ is not injective. Then, there are two points $\af,\be\in M_{s''}$ satisfying $\widetilde{X}_{s''}(\af)=\widetilde{X}_{s''}(\be)$. The convergence of $\{\widetilde X_{s_n}\}_{n\in\natl}$ to $\widetilde X_{s''}$, uniformly in compact subsets of $\real^3$, together with the interior maximum principle, assures that there exist open neighbourhoods $N(\af)$ and $N(\be)$ of $\af$ and $\be$, respectively, such that $\widetilde X_{s''}(N(\af))=\widetilde X_{s''}(N(\be))$. So, the {\it image set} $\widetilde X_{s''}(\widetilde M_{s''})$ is an embedded minimal surface (with boundary) and $\widetilde X_{s''}:\widetilde M_{s''}\longrightarrow\widetilde X_{s''}(\widetilde M_{s''}) $ is a finitely sheeted covering map. As $\widetilde X_{s''}$ is one-to-one in a neighbourhood of the ends, then we deduce that $ \widetilde X_{s''}$ is injective. This contradiction shows that $\cal I$ is closed.

Thus, an elementary connectedness argument gives ${\cal I}=[0,1)$, which implies that the piece of the SC surface contained in the upper half-space is embedded. Using the symmetry about the $x_1$-axis, we conclude the embeddedness of the whole surface.

\hfill q.e.d.
\ \\
\ \\
\ \\
{\bf References}
\ \\
\begin{description}
\itemsep = 0.0 pc
\parsep  = 0.0 pc
\parskip = 0.0 pc
\item{[B]} Bredon, G.E.: Topology and Geometry. Grad. Text in Math. Springer Verlag, New York {\bf 139} (1993)
\item{[C]} Costa, C.: Example of a complete minimal immersion in $\real^3$ of genus one and three embedded ends. Bol. Soc. Bras. Mat. {\bf 15} 41--54 (1984)
\item{[Ch]} Cheng, S.Y.: Eigenfunctions and Nodal Sets. Comment. Math. Helvetici {\bf 51} 43--55 (1976)
\item{[FM]} Ferrer, L., Mart\ih n, F.: Minimal surfaces with helicoidal ends. Preprint math.DG/0405588 (2004).
\item{[HK]} Hoffman, D., Karcher, H.: Complete embedded minimal surfaces of finite total curvature. Encyclopedia of Math. Sci., Springer Verlag {\bf 90} 5--93 (1997) 
\item{[HKW]} Hoffman, D., Karcher, H., Wei, F.: The singly periodic genus-one helicoid. Comment. Math. Helv. {\bf 74} 248--279 (1999)
\item{[HM]} Hoffman, D., Meeks, W.: Embedded minimal surfaces of finite topology. Ann. Math. {\bf 131} 1--34 (1990)  
\item{[K]} Karcher, H.: Construction of minimal surfaces. Surveys in Geometry, University of Tokyo 1--96 (1989) and Lecture Notes {\bf 12}, SFB256, Bonn (1989) 
\item{[LM]} L\'opez, F. J., Mart\ih n, F.: Complete minimal surfaces in $real^3$. Publicacions Matem\`ati- ques {\bf 43} 341--449 (1999)
\item{[MR]} Meeks, W. H., Rosenberg, H.: The global theory of doubly periodic minimal surfaces. Invent. Math. {\bf 97} 351--379 (1989)
\item{[N]} Nitsche, J. C. C.: Lectures on minimal surfaces. Cambridge University Press, Cambridge, (1989)
\item{[O]} Osserman, R.: A survey of minimal surfaces. Dover, New York, 2nd ed (1986)
\item{[R]} Riemann, B.: \"Uber die Fl\"ache vom Kleinsten Inhalt bei gegebener Begrenzung. Abh. K\"onigl. Ges. d. Wiss. G\"ottingen, Mathem. Cl. {\bf 13} 3--52 (1867)
\item{[RB1]} Ramos Batista, V.: Construction of new complete minimal surfaces in $\real^3$ based on the Costa surface. Doctoral thesis, University of Bonn, Germany (2000)
\item{[RB2]} Ramos Batista, V.: Accurate evaluation of elliptic integrals. Technical Report 61/03, University of Campinas (2003) 
\\
http://www.ime.unicamp.br/rel\_pesq/2003/rp61-03.html
\item{[RB3]} Ramos Batista, V.: Polynomial estimates in the unitary interval. Technical Report 60/03, University of Campinas (2003) 
\\
http://www.ime.unicamp.br/rel\_pesq/2003/rp60-03.html
\item{[Sk]} Scherk, H. F.: Bemerkungen \"uber die kleinste Fl\"ache innerhalb gegebener Grenzen. Journal f\"ur die reine und angewandte Mathematik {\bf 13} 185--208 (1835) 
\item{[Sn]} Schoen, R.: Uniqueness, symmetry and embeddedness of minimal surfaces. J. Differential Geom. {\bf 18} 791--809 (1983)
\item{[Sz]} Schwarz, H. A.: Gesammelte Mathematische Abhandlungen. Springer, Berlin (1890) 
\item{[T1]} Traizet, M.: Construction de surfaces minimales en recollant des surfaces de Scherk. Ann. Inst. Fourier, Grenoble, {\bf 46} 1385-1442 (1996)
\item{[T2]} Traizet, M.: An embedded minimal surface with no symmetries. J. Differential Geom. {\bf 60} 103--153 (2002)
\end{description}
\end{document}